\documentclass{article}

\usepackage{hyperref}
\usepackage{graphicx,epsfig}
\usepackage{amsmath,amssymb}
\usepackage{placeins}
\usepackage{algorithmic}
\usepackage{algorithm}
\usepackage[scriptsize]{subfigure}
\usepackage{url}
\usepackage{listings}
\usepackage{verbatim}
\usepackage{colortbl}
\usepackage{subfigmat}
\usepackage{lineno,hyperref}
\usepackage{numcompress}

\modulolinenumbers[5]

\def\R{\mathbb R}
\def\G{\Omega}

\newcommand{\ms}{\medskip}
\newcommand{\tabincell}[2]{\begin{tabular}{@{}#1@{}}#2\end{tabular}}

\newtheorem{Property}{\bf Property}

\begin{document}

\title{A hybrid recursive multilevel incomplete factorization preconditioner for solving  general linear systems}

\author{Yiming Bu\footnotemark[1]~\footnotemark[2]
\and
Bruno Carpentieri\footnotemark[1] \and
Zhaoli Shen\footnotemark[1]~\footnotemark[2]
\and
Tingzhu Huang\footnotemark[2]
}

\renewcommand{\thefootnote}{\fnsymbol{footnote}}

\footnotetext[1]{Johann Bernoulli Institute for Mathematics and Computer Science - University of Groningen, 9747 AG Groningen, The Netherlands.}
\footnotetext[2]{School of Mathematical Sciences, University of Electronic Science and Technology of China, Chengdu, Sichuan 611731, China.}

\maketitle

\makeatother


\begin{abstract}
In this paper we introduce an algebraic recursive multilevel incomplete factorization preconditioner, based on a distributed Schur complement formulation, for solving general linear systems. The novelty of the proposed method is to combine factorization techniques of both implicit and explicit type, recursive combinatorial algorithms, multilevel mechanisms and overlapping strategies to maximize sparsity in the inverse factors and consequently reduce the factorization costs. Numerical experiments demonstrate the good potential of the proposed solver to precondition effectively general linear systems, also against other state-of-the-art iterative solvers of both implicit and explicit form.
\end{abstract}

\ms 


\section{Introduction}\label{sec:1}

\label{sec:one}

Krylov subspace methods may be considered the method of choice for solving large and sparse systems of linear equations arising from the discretization of (systems of) partial differential equations on modern parallel computers. This class of algorithms are iterative in nature. 
At every step $k$, they compute the approximate solution $x_k$ of a linear system $Ax=b$
from the Krylov subspace of dimension~$k$
\[
	K_k(A,b)=span\large\{ b, Ab, A^2b, \ldots, A^{k-1}b\},
\]
according to different criteria for each given method.
The computation requires matrix-vector products with the coefficient matrix 
$A$ plus vector operations, thus potentially reducing the cumbersome costs of sparse direct solvers on large problems, especially in terms of memory. 
All of the iterative Krylov methods converge rapidly if $A$ is somehow close to the identity. Therefore, it is common replacing the original system $Ax=b$ by
\begin{equation}\label{eq:2}
M^{-1}Ax=M^{-1}b,
\end{equation}
or
\begin{equation}\label{eq:3}
AM^{-1}y=b,~~x=M^{-1}y,
\end{equation}
for a nonsingular matrix $M\approx A$. Systems~(\ref{eq:2}) and~(\ref{eq:3}) 
are referred to as {\it{left}} and {\it{right preconditioned}} systems, respectively, 
and $M$ as the {\it{preconditioner matrix}}.
In the case $M$ is factorized as the product of two sparse matrices, $M=M_1M_2$,
like in the Hermitian and positive definite case, one might solve the modified 
linear system
\begin{equation}\label{eq:4}
M_1^{-1}AM_2^{-T}y=M_1^{-1}b,~~x=M_2^{-T}y.
\end{equation}
If one may choose $M$ so that $M^{-1}A$, $AM^{-1}$ or
$M_1^{-1}AM_2^{-T}$ approximate the identity, and linear systems with $M$ or
with $M_1$ and $M_2$ as coefficient matrices are easy to invert, it is more efficient
to apply a Krylov subspace method to the modified linear system. 

Optimal analytic preconditioners based on low order discretizations, nearby equations
that are simple to solve, or similar ideas have been proposed in the literature for 
specific problems. However, the problem-specific approach is generally sensitive to
the characteristics of the underlying operator and to the details of the
geometry. 
In this study, we pursue an algebraic approach where the preconditioner $M$ is computed
only from the coefficient matrix $A$. Although not optimal for any specific
problem, algebraic methods are universally applicable, they can be adapted
to different operators and to changes in the geometry by tuning a few parameters,
and are well suited for solving irregular problems defined on unstructured grids.

Roughly speaking, most of the existing techniques can be divided into either
implicit or explicit form. A preconditioner of {\em implicit} form is defined
by any nonsingular matrix $M\approx A$, and requires to solve an extra linear system with $M$ at each
step of an iterative method. The most important example in this class is  
represented by the Incomplete $LU$ decomposition (ILU), where $M$ is implicitly defined as $M=\bar L \bar U$, $\bar L$ and $\bar U$ being triangular matrices that approximate the exact $L$ and $U$ factors of $A$ according to a prescribed dropping strategy adopted during the Gaussian elimination process. These methods are considered amongst the most reliable
in a general setting. Well known theoretical results on the existence and the stability of
the factorization can be proved for the class of $M$-matrices~\cite{mevd:77}, and recent
studies are involving more general matrices, both structured and unstructured.
The quality of the factorization on difficult problems can be enhanced by using several techniques such as reordering, scaling, diagonal shifting, pivoting and condition estimators~(see e.g.~\cite{duko:99a,ARMS,mmbw:00,BolS06,cabo:12}). As a result of this active development,
in the last years successful results are reported with ILU-type preconditioners in many areas that were of exclusive domain of direct solution methods like in circuits simulation, power
system networks, chemical engineering plants modelling, graphs and other problems
not governed by PDEs, or in areas where direct methods have been traditionally preferred,
like in structural analysis, semiconductor device modelling and computational fluid
dynamics applications~(see e.g.~\cite{saad:2005,boll:2003,Benzi02,Manguoglu2011,SaSoTo02}). 
One problem with ILU-techniques is the severe degradation of performance 
observed on vector, parallel and GPUs machines, mainly due to the sparse triangular solves~\cite{lisa:13}. 
In some cases, reordering techniques may help to introduce nontrivial parallelism. However, parallel orderings may sometimes degrade the convergence rate, while more fill-in diminishes the overall parallelism of the solver~\cite{dume:89}.

{\em Explicit} preconditioning tries to mitigate such difficulties by approximating directly $A^{-1}$, as the product $M$ of sparse matrices, so that the
preconditioning operation  reduces to forming one (or more) sparse matrix-vector product,
and consequently the application of the preconditioner may be easier to parallelize and numerically stable.
Some methods can also  perform the construction phase in parallel~\cite{huckle10efficient,cdgs:05,BFSAI,pash:14,pash:14a}; additionally,
on certain indefinite problems with large nonsymmetric parts, the explicit approach 
can provide better results than ILU techniques (see e.g.~\cite{chsa:98,carp:07d,huckle10smoothing}).
In practice, however, some questions need to be addressed.
The computed matrix $M$ could be singular, and the construction cost is typically 
much higher than for ILU-type methods, especially for sequential runs.
The main issue is the selection of the non-zero pattern of $M$. The idea is to keep $M$
reasonably sparse while trying to capture the `large' entries of the inverse, which are
expected to contribute the most to the quality of the preconditioner.
On general problems it is difficult to determine the best structure for $M$ in advance,
and the computational and storage costs required to achieve the same rate of convergence of preconditioners given in implicit form may be prohibitive in practice. 

In this study, we present an algebraic multilevel solver for preconditioning general nonsymmetric linear systems which attempts to combine characteristics of both approaches. 
Assuming that the matrix $A$ admits the factorization $A=LU$, with $L$ a unit lower and $U$ an upper triangular matrix, our method approximates the inverse factors $L^{-1}$ and $U^{-1}$. Sparsity in the approximate inverse factors is maximized by employing recursive combinatorial algorithms. Robustness is enhanced by combining the factorization with recently developed overlapping strategies and by using efficient local solvers.

The paper is organized as follows. In Section~\ref{sec:3} we describe the proposed multilevel preconditioner. In Section~\ref{sec:4} we show how to combine our preconditioner with overlapping strategies, and in Section~\ref{sec:5} we assess its overall performance by showing several numerical experiments on realistic matrix problems, also against other state-of-the-art solvers. Finally, in Section~\ref{sec:6} we conclude the study with some remarks and perspectives for future work.

\section{The AMES solver}\label{sec:3}

Let 
\begin{equation}\label{eq:1}
	Ax=b
\end{equation}
be a $n \times n$ general linear system with nonsingular, possibly indefinite and nonsymmetric matrix $A = \{a_{ij}\} \in \R^{n \times n}$, and vectors $x, b \in \R^{n}$. 
We assume that $A$ admits for a triangular decomposition
\[
	A=LU
\]
and we precondition system~(\ref{eq:1}) as
\[
	M_L A M_U y = M_L b
\]
with $M_L \approx L^{-1}$ and $M_U  \approx U^{-1}$, clearly preserving symmetry and/or positive definiteness of $A$. This approach of preconditioning linear systems has been extensively investigated in a series of papers by Kolotilina and Yeremin ~\cite{koye:93,koye:95,koyn:00,koyn:99}, who prescribed the nonzero pattern of the inverse factors $M_L$ and $M_U$ of $A$ in advance equal to the pattern of the lower and upper triangular part of $A+A^T$, respectively, and determined the entries of $M_L$ and $M_U$ explicitly by solving linear equations involving the principal submatrices of $A$ ({\it{the `FSAI' preconditioner}}). Chow suggested to use as pattern for the inverse factors the structure of the lower and upper triangular part of $(A+A^T)^p$, where $p$ is a positive integer~\cite{chow:00a,chow:01,yuni:04}. The larger $p$, in general the higher the quality of the computed preconditioner, although the construction, storage and application costs tend to increase rapidly with $p$. Blocking and adaptive strategies have been recently studied to overcome these problems~\cite{BFSAI,gBFSAI,JANNA}. Benzi and T\accent23uma proposed to compute the entries of matrices $M_L$ and $M_U$ by means of a (two-sided) Gram-Schmidt orthogonalization process with respect to the bilinear form associated with $A$, and to determine the best structure for the inverse factors dynamically, during the construction ({\it{the `AINV' preconditioner}}). Sparsity is preserved in the process by discarding elements having magnitude smaller than a given positive threshold~\cite{bemt:96,betu:98b}. 

In this study we analyse multilevel mechanisms, recursive combinatorial algorithms and overlapping techniques, combined with efficient local solvers, to enhance robustness and reduce costs for the approximation of the inverse factors. We refer to the resulting preconditioner as AMES~(Algebraic Multilevel Explicit Solver). It is easier to describe the AMES method by using graph notation, dividing the solution of system~(\ref{eq:1}) in five distinct phases:

\begin{enumerate}

  \item

a \emph{scale phase}, where the coefficient matrix $A$ is scaled by rows and columns so that the largest entry of the scaled matrix has magnitude smaller than one;

  \item

a \emph{preorder phase}, where the structure of $A$ is used to compute a suitable ordering that maximizes sparsity in the approximate inverse factors;

  \item

an \emph{analysis phase}, where the sparsity preserving ordering is analyzed and an efficient data structure is generated for the factorization;

  \item

a \emph{factorization phase}, where the nonzero entries of the preconditioner are computed;

  \item

a \emph{solve phase}, where all the data structures are accessed for solving the linear system.

\end{enumerate}

Below we describe each phase separately.

\subsection{Scale phase.}
We initially scale system~(\ref{eq:1}) by rows and columns as
\begin{equation}\label{eq:linscal}
    D_1^{1/2}Ay=D_1^{1/2}b,~~y=D_2^{1/2}x,
\end{equation}
where the ${n \times n}$ diagonal scaling matrices $D_1$ and $D_2$ have the form
\[
D_1 (i,j) = \left\{ {\begin{array}{*{20}c}
   {\frac{1}{{\mathop {\max \left| {a_{ij} } \right|}\limits_{i} }}} & , \rm{~~if~i~=~j} \\
      & \\
   0 & , \rm{~~if~i~\ne~j} \\
\end{array}} \right.,
~~~D_2 (i,j) = \left\{ {\begin{array}{*{20}c}
   {\frac{1}{{\mathop {\max \left| {a_{ij} } \right|}\limits_{j} }}} & , \rm{~~if~i~=~j} \\
      & \\
   0 & , \rm{~~if~i~\ne~j} \\
\end{array}} \right.~.
\]
For simplicity, we still refer in this paper to the scaled system~(\ref{eq:linscal}) as $Ax=b$.

\subsection{Preorder phase.}
We use standard notation of graph theory to describe this computational step. We denote by $\G(\tilde A)$ the undirected graph associated with the matrix
\[
\tilde A = \left\{ \begin{array}{l}
 A ,~~~~~~~~{\rm{~if~}} A {{\rm{~is~symmetric,}}} \\
 A + A^T , {\rm{~if~}} A {{\rm{~is~nonsymmetric.}}} \\
 \end{array} \right.
\]
First, $\G(\tilde A)$ is partitioned into $p$ non-overlapping subgraphs $\G_i$ of roughly equal size by using the multilevel graph partitioning algorithms available in the Metis package~\cite{METIS}.
For each partition $\G_i$ we distinguish two disjoint sets of nodes (or vertices):~{\it{interior nodes}} that are connected only to nodes in the same partition, and {\it{interface nodes}} that straddle between two different partitions; the set of interior nodes of $\G_i$ form a so called {\it{separable}} or {\it{independent cluster}}. Upon renumbering the vertices of $\G$ one cluster after another, followed by the interface nodes as last, and permuting $A$ according to this new ordering, a block bordered linear system is obtained, with coefficient matrix of the form
\begin{equation}\label{eq:perm}
\tilde{A}=P^T AP = \left( {\begin{array}{*{20}c}
   B & F  \\
   E & C  \\
 \end{array} } \right) =
\left( {\begin{array}{*{20}c}
   {B_1 } & {} & {} & {} & {F_1 }  \\
   {} & {B_2 } & {} & {} & {F_2 }  \\
   {} & {} &  \ddots  & {} &  \vdots   \\
   {} & {} & {} & {B_p } & {F_p }  \\
   {E_1 } & {E_2 } &  \cdots  & {E_p } & {C}  \\
\end{array}} \right).
\end{equation}
In~\eqref{eq:perm}, each diagonal block $B_i$ corresponds to the interior nodes of $\G_i$, and the blocks $E_i$ and $F_i$ correspond to the interface nodes of $\G_i$; the block $C$ is associated to the mutual interactions between the interface nodes. In our multilevel scheme we apply the same block downward arrow structure to the diagonal blocks $B_i$ of $\tilde{A}$;
the procedure is repeated recursively until a maximum number of levels is reached, or until the blocks at the last level are sufficiently small to be easily factorized. As an example, in Figure~\ref{fig:permuted}(b) we show the structure of the sparse matrix \emph{rdb2048} from Tim Davis’ matrix collection~\cite{timdavis} after three reordering levels.

To reduce factorization costs, a similar permutation is applied to the Schur complement matrix $S=C-EB^{-1}F$ as follows
\begin{equation}\label{eq:permAS}
\tilde{S}= \left( {\begin{array}{*{20}c}
   {B_{S1} } & {} & {} & {} & {F_{S1} }  \\
   {} & {B_{S2} } & {} & {} & {F_{S2} }  \\
   {} & {} &  \ddots  & {} &  \vdots   \\
   {} & {} & {} & {B_{Sp} } & {F_{Sp} }  \\
   {E_{S1} } & {E_{S2} } &  \cdots  & {E_{Sp} } & {C_S}  \\
\end{array}} \right).
\end{equation}

\subsection{Analysis phase.}

In the analysis phase, a suitable data structure for storing the linear system is defined, allocated and initialized. We use a tree structure to store the block bordered form~(\ref{eq:perm}) of $\tilde{A}$. The root is the whole graph $\G$, and the leaves at each level are the independent clusters of each subgraph. Each node of the tree corresponds to one partition $\G_i$ of $\G(\tilde A)$, or equivalently to one block $B_i$ of matrix $\tilde{A}$. The information stored at each node are the entries of the off-diagonal blocks $E$ and $F$ of $B_i's$ father, and those of the block $C$ of $B_i$ after its permutation, except at the last level of the tree where we store the entire block $B_i$. All these matrices are represented in the computer memory using a compressed sparse row storage format, except for blocks $F_i$ that are stored in compressed sparse column format. Blocks $E_i$ and $F_i$ can be very sparse; many of their rows and columns can be zero, and this leads to a significant saving of computation.

\subsection{Factorization phase.}\label{sec:Factorization}
The approximate inverse factors $\tilde{L}^{-1}$ and $\tilde{U}^{-1}$ of~$\tilde{A}$ write in the following form
\begin{equation}\label{eq:nbyninv}
\tilde{L}^{-1} \approx
\left(
{\begin{array}{*{20}c}
   {U_1^{ - 1} } & {} & {} & {} & {W_1 }  \\
   {} & {U_2^{ - 1} } & {} & {} & {W_2 }  \\
   {} & {} &  \ddots  & {} &  \vdots   \\
   {} & {} & {} & {U_p^{ - 1} } & {W_p }  \\
   {} & {} & {} & {} & {U_S^{ - 1} }  \\
\end{array}} \right) , ~~
\tilde{U}^{-1}  \approx
\left( {\begin{array}{*{20}c}
   {L_1^{ - 1} } & {} & {} & {} & {}  \\
   {} & {L_2^{ - 1} } & {} & {} & {}  \\
   {} & {} &  \ddots  & {} & {}  \\
   {} & {} & {} & {L_p^{ - 1} } & {}  \\
   {G_1 } & {G_2 } &  \cdots  & {G_p } & {L_S^{ - 1} }  \\
\end{array}} \right)
\end{equation}
where
\begin{equation}\label{eq:borders}
B_i  = L_i U_i, W_i  =  - U_i^{ - 1} L_i^{ - 1} F_i U_S^{ - 1}, G_i
=  - L_S^{ - 1} E_i U_i^{ - 1} L_i^{ - 1}
\end{equation}
and $L_S, U_S$ are the triangular factors of the Schur complement matrix
\begin{equation}\label{eq:schurbuild}
S = C  - \sum\limits_{i = 1}^p {E_i B_i^{ - 1} F_i }.
\end{equation}

\begin{figure}
\subfigure[The original structure of the \emph{rdb2048} matrix.]{\includegraphics[width=3.5cm]{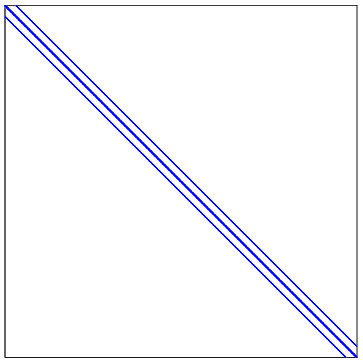}}
\hfill
\subfigure[The structure of \emph{rdb2048} after permutation.]{\includegraphics[width=3.5cm]{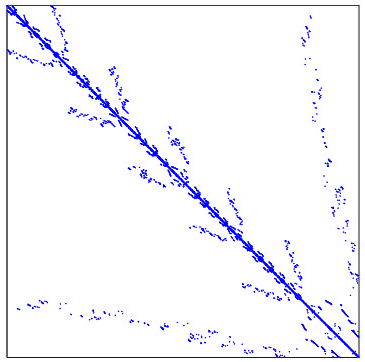}}
\hfill
\subfigure[The structure of the inverse factor.  In red are displayed the entries actually stored.]{\includegraphics[width=3.5cm]{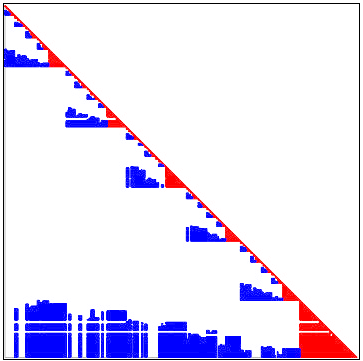}}
\caption{Structure of the multilevel inverse-based factorization for the matrix \emph{rdb2048}.}~\label{fig:permuted}
\end{figure}

Some fill-in may occur in $\tilde{L}^{-1}$ and $\tilde{U}^{-1}$ during the factorization, but only within the nonzero blocks. This two-level reordering scheme was used in the context of factorised approximate inverse methods for the parallelization of the AINV preconditioner in~\cite{bemt:07}. Differently from~\cite{bemt:07}, we apply the arrow structure~(\ref{eq:perm}) recursively to the diagonal blocks and to the first level Schur complement as well, to gain additional sparsity. The multilevel factorization algorithm requires to invert only the last level blocks and the small Schur complements at each reordering level; the blocks $W_i$, $G_i$ do not need to be assembled explicitly, as they may be applied using Eqn~(\ref{eq:borders}). For the \emph{rdb2048} problem, in Figure~\ref{fig:permuted}(c) we display in red the actual extra storage required by the exact multilevel inverse factorization in addition to matrix $A$; these represent only $34 \%$ of the total nonzeros of $A$. From the knowledge of the red entries, the blue ones can be retrieved from Eqn~(\ref{eq:borders}), using the off-diagonal blocks of $A$.
We also permute the large Schur complement at the first level into a block bordered structure, until we reach a maximal number of levels or a given minimal size. The last-level matrix is inverted inexactly. An inexact solver is also used to factorize the last-level blocks $B_i$ in~(\ref{eq:schurbuild}).

\subsection{Solve phase.}\label{sec:Solve}
In the solve phase, the multilevel factorization is applied at every iteration step of a Krylov method for solving the linear system. Notice that the inverse factorization of~$\tilde{A}$ may be written as
\begin{equation}\label{eq:2by2inv}
\left(P A P^T\right)^{-1}  =
\left( {\begin{array}{*{20}c}
   { U^{ - 1} } & W  \\
   0 & {U_S^{ - 1} }  \\
\end{array}} \right) \times \left( {\begin{array}{*{20}c}
   { L^{ - 1} } & 0  \\
   G & {L_S^{ - 1} }  \\
\end{array}} \right)
\end{equation}
where $W =  - U^{ - 1}  L^{ - 1} F U_S^{ - 1},~G =  - L_S^{ - 1} E U^{ - 1} L^{ - 1}$, and $L_S,~U_S$ are the inverse factors of the Schur complement matrix $S  = C  - E B^{ - 1} F$.

From Eqn.~(\ref{eq:2by2inv}), we obtain the following expression for the exact inverse
\begin{equation}\label{eq:exinv}
\left( {\begin{array}{*{20}c}
   {B^{ - 1}  + B^{ - 1} FS^{ - 1} EB^{ - 1} } & { - B^{ - 1} FS^{ - 1} }  \\
   { - S^{ - 1} EB^{ - 1} } & {S^{ - 1} }  \\
\end{array}} \right).
\end{equation}
We can derive preconditioners from Eqn.~(\ref{eq:exinv}) by computing approximate solvers $\tilde B^{ - 1}$ for $B$ and $\tilde S^{ - 1}$ for $S$. Hence the preconditioner matrix $M$ will have the form 
\[
M = \left( {\begin{array}{*{20}c}
   {\tilde B^{ - 1}  + \tilde B^{ - 1} F \tilde S^{ - 1} E \tilde B^{ - 1} } & { - \tilde B^{ - 1} F \tilde S^{ - 1} }  \\
   { - \tilde S^{ - 1} E \tilde B^{ - 1} } & {\tilde S^{ - 1} }  \\
\end{array}} \right).
\]
and the preconditioning operation
$\left[ {\begin{array}{*{20}c}
   {y_1 }  \\
   {y_2 }  \\
\end{array}} \right] = M\left[ {\begin{array}{*{20}c}
   {x_1 }  \\
   {x_2 }  \\
\end{array}} \right]$
writes as Algorithm~\ref{alg:prec}. Notice that Algorithm~\ref{alg:prec} is called recursively at lines~1-3, as $\tilde B$ and $\tilde S$ also have a block bordered structure upon permutation.
\begin{algorithm}[!ht]
\caption{\label{alg:prec}{\it The preconditioning operation in the AMES solver}.}
\begin{algorithmic}[1]
\STATE $p_1=\tilde B^{-1} x_1$
\STATE $[p_2,p_3]=\tilde S^{-1}[E \cdot p_1, x_2]$
\STATE $[p_4,p_5]=\tilde B^{-1}[F \cdot p_2, F \cdot p_3]$
\STATE $y_1=p_1+p_4-p_5$
\STATE $y_2=p_3-p_2$
\end{algorithmic}
\end{algorithm}

~

\section{Combining the AMES solver with overlapping}\label{sec:4}

In~\cite{overlapping}, Grigori, Nataf and Qu have introduced an~\emph{overlapping technique} to enhance the robustness of multilevel incomplete LU factorization preconditioning computed from matrices reordered in arrow form, e.g. using the nested dissection method by George~\cite{geli:81}. The multilevel mechanism incorporated in the AMES preconditioner described in the previous section is based on a nested dissection-like ordering, and thus it can easily accomodate for overlapping. We have tested this idea in our numerical experiments, and in this section we shortly describe the procedure adopted. The results of our experiments are reported in Section~\ref{sec:5}.

\ms

{\subsection{Background}}
Let $\Omega=(V(\Omega),E(\Omega))$ be the graph of $A$, $V(\Omega)$ denoting the set of vertices and $E(\Omega)$ the set of edges in $\Omega$.
If the graph is directed, we denote an edge of $E$ issuing from vertex $u$ to vertex $v$ as $(u,v)$; $u$ is called a predecessor of $v$, and $v$ a successor of $u$. 
If the graph is undirected, we denote the edges of $E$ by non-ordered pairs $\{u,v\}$;  $u$ is called a neighbour of $v$. As in the previous section, we assume that $\G$ is partitioned into $p$ independent non-overlapping subgraphs $\Omega_{1}$, $\ldots$, $\Omega_{p}$, and we call $S$ the set of separator nodes, straddling between two different partitions. Goal of overlapping is to extend each independent set of $\Omega$ by including its direct neighbours, similarly to the overlapping idea used in other domain decomposition methods, for example in the restricted additive Schwarz method~\cite{quva:99,SAAD-BOOK}.

Following~\cite{overlapping}, we denote by $V(\Omega_{i,ext})$ the separator
nodes that are successors of $\Omega_i$,
\begin{equation}\label{iext}
V(\Omega_{i,ext})=\{v\in V(S)|\exists u\in V(\Omega_i), (u,v)\in E(\Omega)\} \subset V(S),
\end{equation}
and by $V(\Omega_{ext})$ the complete set of successor nodes of all the subdomains
\begin{equation}\label{ext}
V(\Omega_{ext})=\bigcup_{i=1:p}V(\Omega_{i,ext}).
\end{equation}
Then $\Omega_i$ is extended to the set $\hat{\Omega}_i$ as
\begin{equation}\label{hati}
V(\hat{\Omega}_i)=V(\Omega_i)\cup V(\Omega_{i,ext}), ~ i = 1,\ldots,p,
\end{equation}
and the separator $S$ is extended to $\hat S$ by adding the successors of nodes in $V(\Omega_{ext})$, that is
\begin{equation}\label{hats}
V(\hat{S})=V(S)\cup \{v\in V(\Omega_i), i=1,\ldots,p ~|~ \exists u\in V(\Omega_{ext}), (u,v)\in E(\Omega)\}.
\end{equation}

Using this notation, the overlapped graph of $A$, $\tilde \Omega = (V(\tilde \Omega), E(\tilde \Omega))$, is introduced as follows. First define the overlapped subgraph $\tilde \Omega_i$ and the overlapped separator $\tilde S$ as, respectively, 
$$V(\tilde{\Omega}_i)=\{(x,i):x\in\hat{\Omega}_i\},$$
$$V(\tilde{S})=\{(x,s):x\in\hat{S}\}.$$
For simplicity we refer to $(x,i)$ as $x_i$.
Then the vertex set $V(\tilde \Omega)$ of the overlapped graph $\tilde{\Omega}$ is formed by the disjoint union of the $V(\hat{\Omega_i})$'s and of $V(\hat{S})$ as
\begin{equation}\label{tilde}
V(\tilde{\Omega})= \left( \bigcup_{i\in 1:p} V(\hat \Omega_i) \right) \cup V(\hat S).
\end{equation}
Recall that, given the union $B$ of a family of sets indexed by the index set $I$,
$$B=\bigcup_{i\in I}A_i=\bigcup_{i\in I}\{x:x\in A_i\},$$
their disjoint union $C$ is defined as the set
$$C=\bigcup_{i\in I}\{(x,i):x\in A_i\}.$$
At this stage, it is useful to introduce the two projection operators $\Pi_{1}$ and $\Pi_{2}$ such that
$$\Pi_{1}: (x,i)\mapsto x$$
and
$$\Pi_{2}: (x,i)\mapsto i.$$
With this notation, the set of edges of the overlapped graph $\tilde{\Omega}$
is defined according to their projection onto the original graph as follows
\begin{equation}\label{edgei}
E(\tilde{\Omega}_{i})=\{(u_i,v_i)|u_i\in V(\tilde{\Omega}_{i}),v_i\in V(\tilde{\Omega}_{i}), (\Pi_1 (u_i),\Pi_1 (v_i))\in E(\Omega)\},
\end{equation}
\begin{equation}\label{edges}
E(\tilde{S})=\{(u_s,v_s)|u_s\in V(\tilde{S}),v_s\in V(\tilde{S}), (\Pi_1 (u_s),\Pi_1 (v_s))\in E(\Omega)\},
\end{equation}

$$E(\tilde{\Omega}_{i},\tilde{S})=\{(u_i,v_s)|u_i\in V(\tilde{\Omega}_{i}),v_s\in V(\tilde{S}), (\Pi_1 (u_i),\Pi_1 (v_s))\in E(\Omega),$$
\begin{equation}\label{edgeis}
\nexists v_i\in V(\tilde{\Omega}_{i}), (u_i,v_i)\in E(\tilde{\Omega}_{i})\},
\end{equation}

$$E(\tilde{S},\tilde{\Omega}_{i})=\{(u_s,v_i)|u_s\in V(\tilde{S}),v_i\in V(\tilde{\Omega}_{i}), (\Pi_1 (u_s),\Pi_1 (v_i))\in E(\Omega),$$
\begin{equation}\label{edgesi}
\nexists v_s\in V(\tilde{S}), (u_s,v_s)\in E(\tilde S)\}.
\end{equation}

The following property, established in~\cite{overlapping}, ensures an equivalence between the equations of the overlapped system and those of the original system. 

\begin{Property}
Let $\Omega$ be the associated directed graph of a given system of linear equations and $u$ be a vertex of $V(\Omega)$. Let $\tilde{\Omega}$ be the overlapped graph, and let $u_i$ be a vertex of $V(\tilde{\Omega})$ such that $\Pi_{1}(u_i)=u \in V(\Omega)$. For each edge $(u,v)\in E(\Omega)$, there is a unique $v_j\in V(\tilde{\Omega})$ such that we have both $\Pi_{1}(v_j)=v$ and $(u_i,v_j)\in E(\tilde{\Omega})$.
\end{Property}

This property shows that there exists a bijection between the nonzeros of the equation corresponding to vertex $u$ in the original system and the nonzeros of the equation corresponding to its dual $u_i$, where $\Pi_{1}(u_i)=u$. From a matrix viewpoint, to each nonzero entry $\tilde{a}_{u_i,v_i}$ in the overlapped matrix there is a unique nonzero entry $a_{u,v}$ in the original matrix, where $\Pi_{1}(u_i)=u$ and $\Pi_{1}(v_i)=v$. Therefore there is a one-to-one correspondence between equations in the original system and those in the overlapped system. By solving the overlapped system, we can automatically obtain the solution of the original system.

\ms

{\subsection{Example}} 
We display a simple example from~\cite{overlapping} to describe shortly how the overlapping procedure works in practice. We consider a $5\times 5$ matrix having the structure shown in Figure~\ref{fig:smalloverlapping}(a). The graph consists of two independent subgraphs $\Omega_1=\{1,2\}$, $\Omega_2=\{3\}$ and one separator $S=\{4,5\}$. We just pick the first subgraph and the separator set to explain. Separator nodes that are successors of $\Omega_1$ are the set
\[V(\Omega_{1,ext})=\{4,5\}\]
and we have
\[ V(\hat{\Omega}_1)=V(\Omega)_1\cup V(\Omega_{1,ext})=\{1,2,4,5\},\]
so that
\[V(\tilde{\Omega}_1)=\{1_1,2_1,4_1,5_1\}.\]
Analogously, \[V(\Omega_{2,ext})=\{4,5\}\] and
\[V(\Omega_{ext})=\Omega_{1,ext}\cup \Omega_{2,ext}=\{4,5\}.\]
Next, we compute the overlapped separator set $\tilde S$. The vertices of $V(\Omega_1)$ and $V(\Omega_2)$ directed by $V(\Omega_{ext})$ are $\{1,3\}$, so
\[V(\hat{S})=V(S)\cup\{1,3\}=\{4,5,1,3\}\]
and
\[V(\tilde{S})=\{4_s,5_s,1_s,3_s\}.\]

According to Eqns.~(\ref{edgei})-(\ref{edgesi}), the edges of the overlapped subdomain $E(\tilde{\Omega}_1)$ are defined based on their projection onto the original graph. The first diagonal block of the overlapped matrix is formed by picking the $V(\hat{\Omega}_1)= \{1,2,4,5\}$ rows and columns of the original matrix
$$
\bordermatrix{%
& 1 & 2 & 4 & 5 \cr
1&\diamond & \diamond & \diamond & \diamond \cr
2&      \diamond & \diamond & \  & \  \cr
4&      \diamond & \  & \diamond & \  \cr
5&      \diamond & \  & \  & \diamond \cr
}.
$$
Similarly for the other two diagonal blocks, and this is shown in  Figure~\ref{fig:smalloverlapping}(b).

From~Eqn.~(\ref{edgeis}), we can construct the edges from $\tilde{\Omega}_1$ to $\tilde{S}$. These are the nonzero entries of the overlapped interface block $\tilde F_1$, adopting the same notation as in (\ref{eq:perm}). We pick the $V(\hat{\Omega}_1)=\{1,2,4,5\}$ rows and $V(\hat S)=\{4,5,1,3\}$ columns of the original matrix, and we set the columns corresponding to the common vertexes of $\hat{\Omega}_1$ and $\hat{S}$ to zeros. In our example this results in zeroing out the columns of $\hat F_1$ indexed by $\hat{\Omega}_1\cup\hat{S}=\{4,5,1\}$, giving
$$\bordermatrix{%
& 4 & 5 & 1 & 3 \cr
1 & \times & \times & \diamond & \  \cr
2 & \  & \  & \diamond & \  \cr
4 & \star & \  & \times &\times \cr
5 & \  & \times & \times &\times \cr
}
\longrightarrow
\bordermatrix{%
& 4 & 5 & 1 & 3 \cr
1 & \  & \  & \  & \  \cr
2 & \  & \  & \  & \  \cr
4 & \  & \  & \  &\times \cr
5 & \  & \  & \  &\times \cr
}.$$
Similar procedure is followed for the other blocks $\tilde F_i$, $\tilde E_i$. Finally, the overlapped matrix has the form given in Figure~\ref{fig:smalloverlapping}(b). The block arrow structure of the original matrix is preserved. However, symmetry is lost and the sparsity pattern also changes significantly.

\begin{figure}
 \begin{center}
\begin{minipage}[t]{10cm}%
\begin{large}
\[
  \left[ {\begin{array}{cc|c|cc}
   a_{11} & a_{12} & ~ & a_{14} & a_{15}\\
   a_{21} & a_{22} & ~ & ~ & ~\\
   \hline
   ~ & ~ & a_{33} & a_{34} & a_{35}\\
   \hline
   a_{41} & ~ & a_{43} & a_{44} & ~ \\
   a_{51} & ~ & a_{53} & ~ & a_{55} \\
  \end{array} } \right]
\]
\end{large}
  \center{(a) The original matrix}
\end{minipage}%

\begin{minipage}[t]{10cm}%
\begin{small}
  \[
  \left[ {\begin{array}{cccc|ccc|cccc}
   a_{11} & a_{12} & a_{14} & a_{15} & ~ & ~ & ~ & ~ & ~ & ~ & ~ \\
  a_{21} & a_{22} & ~ & ~ & ~ & ~ & ~ & ~ & ~ & ~ & ~ \\
   a_{41} & ~ & a_{44} & ~ & ~ & ~ & ~ & ~ & ~ & ~  & a_{43} \\
   a_{51} & ~ & ~ & a_{55} & ~ & ~ & ~ & ~ & ~ & ~  & a_{53} \\
   \hline
    ~ & ~ & ~ & ~ & a_{33} & a_{34} & a_{35} & ~ & ~ & ~ & ~ \\
    ~ & ~ & ~ & ~ & a_{43} & a_{44} & ~ & ~ & ~ & a_{41} & ~ \\
    ~ & ~ & ~ & ~ & a_{53} & ~ & a_{55} & ~ & ~ & a_{51} & ~ \\
    \hline
    ~ & ~ & ~ & ~ & ~ & ~ & ~ & a_{44} & ~ & a_{41} & a_{43} \\
    ~ & ~ & ~ & ~ & ~ & ~ & ~ & ~ & a_{55} & a_{51} & a_{53} \\
    ~ & a_{12} & ~ & ~ & ~ & ~ & ~ & a_{14} & a_{15} & a_{11} & ~ \\
    ~ & ~ & ~ & ~ & ~ & ~ & ~ & a_{34} & a_{35} & ~ & a_{33} \\
  \end{array} } \right]
\]
\end{small}

  \center{(b) The matrix after one-level overlapping}
\end{minipage}
 \end{center}
\caption{Matrix structure before and after applying the overlapping procedure.}\label{fig:smalloverlapping}
\end{figure}

\bigskip

{\subsection{Analysis}}~It is interesting to analyse the effect that overlapping may produce on the AMES method. According to (\ref{hati}) and (\ref{edgei}), the size and the number of nonzeros in each subgraph is increased after overlapping. According to (\ref{edgeis}), the interconnections between subdomains and separator are reduced in the overlapped system, as the original interconnectivities are all removed. The more nodes are added to the subgraphs, the richer they are in terms of information about the system matrix, and a larger performance improvement may be expected. In the overlapped system, each block $\tilde B_i$ has the following structure

$$\bordermatrix{%
& V(\Omega_i) & V(\Omega_{i,ext}) \cr
V(\Omega_i) & B  & E(\Omega_i,\Omega_{i,ext})   \cr
V(\Omega_{i,ext}) & E(\Omega_{i,ext},\Omega_i)  & E(\Omega_{i,ext})  \cr
}.$$

From Eqn.~(\ref{iext}) we see that the set of neighboring nodes $V(\Omega_{i,ext})$  corresponds to the nonzero columns of the block $F_i$, and the nonzero elements of $F_i$
are determined by the set of interconnections $E(\Omega_i,\Omega_{iext})$. Therefore, the more dense and larger the blocks $F_i, i=1:p,$ (that is,  the size of separator) in the original matrix, the more nodes and interconnections are added to subdomains, and a larger reduction of the number of iterations can be achieved.

\bigskip

{\subsection{Algorithmics}} 

The AMES preconditioning algorithm described in Section~\ref{sec:3} with one extra overlapping phase writes as follows:

\ms 

\begin{enumerate}

  \item

a \emph{scale phase}, where the matrix $A$ is scaled by rows and columns so that the largest entry of the scaled matrix has magnitude smaller than one;

  \item

a \emph{preorder phase}, where the structure of $A$ is used to compute a suitable ordering that can maximize sparsity in the approximate inverse factors;

  \item

an \emph{overlap phase}, which extends each block $B_i$ and the Schur complement, and generates the overlapped matrix $\tilde A$ and the right-hand side vector $\tilde b$;

  \item

an \emph{analysis phase}, where the sparsity preserving and overlapping orderings are analyzed and an efficient data structure is generated for the factorization;

  \item

a \emph{factorization phase}, where the entries of $\tilde A$ are processed to explicitly compute the approximate inverse factors;

  \item

a \emph{solve phase}, that accesses all the data structures for solving the overlapped linear system.

  \item
  
a \emph{restriction phase}, that restricts the solution obtained from the overlapped system to the original system, and obtains the solution.
\end{enumerate}

~

\section{Numerical experiments}\label{sec:5}

In this section we present the results of our numerical experiments to illustrate the performance of the AMES preconditioner, also against other state-of-the-art methods and software for solving general linear systems. The selected  matrix problems are extracted from the public-domain matrix repository available at the University of Florida~\cite{timdavis}, and arise from various application fields. We present a summary of the characteristics of each linear system in Table~\ref{tab:problems}.
\begin{table}[!ht]
 \begin{center}
  \begin{tabular}{lclr}
    \hline
      Matrix problem & $n$ & Field & $nnz(A)$  \\
    \hline
     \texttt{orsirr\_1} & ~1,030~  & Oil reservoir simulation  &  6,858  \\
     \texttt{1138\_bus} & ~1,138~  & Bus Power System  &  4,054  \\
     \texttt{bcsstk27} & ~1,224~  & BCS Structural Engineering Matrix   &  28,675  \\
     \texttt{epb0} & ~1,794~  & Plate-fin heat exchanger  &  7,764 \\
     \texttt{cz20468} & 20,468          & Closest Point Method       &  206,076  \\
     \texttt{raefsky3} & 21,200       	& Fluid Structure Interaction &  1,488,768 \\
     \texttt{ABACUS\_shell\_ud} & 23,412   & ABAQUS benchmark &  218,484 \\     \texttt{sme3Db} & 29,067           & 3D structural mechanics problem     &  2,081,063  \\     \texttt{viscoplastic2} & 32,769    & FEM discretization   &  381,326  \\     \texttt{cz40948} & 40,948          & Closest Point Method   &  412,148  \\
     \texttt{rma10} & ~46,835~            & 3D CFD Model           &  2,374,001  \\
     \texttt{finan512} & ~74,752~            & Portfolio optim     &  596,992  \\
     \texttt{helm2d03} & ~392,257~            &  Helmholtz eq. on a unit square     &  2,741,935  \\
     \texttt{parabolic\_fem} & ~525,825~     & Parabolic FEM  &  3,674,625 \\
     
    \hline
  \end{tabular}
  \ms
  \caption{\label{tab:problems} Set and characteristics of the test matrix problems.}
 \end{center}
\end{table}
We applied AMES as a preconditioner for the Generalized Minimal Residual (GMRES) method by Saad and Schultz~\cite{sasc:86}. In all our runs we started the iterative solution from the zero vector and we stopped it when either the initial residual was reduced by twelve orders of magnitude or when no convergence was achieved after 5000 matrix-vector products. To limit  memory costs, we restarted the GMRES algorithm every $500$ iterations. The right-hand side $b$ of the linear system was chosen so that the solution is the vector of all ones, that is $b=Ae$ with $e=[1,...,1]^{T}$. In each run we recorded the following performance measures:

\begin{enumerate}\addtolength{\itemsep}{0.2\baselineskip}
\item the density ratio $\frac{nnz(M_L + M_U)}{nnz(A)}$, that is the ratio between the number of nonzeros in the preconditioner matrix $M=M_UM_L$ versus the number of nonzeros in the coefficient matrix $A$;
\item the number of iterations $Its$ required to reduce the initial residual by 12 orders of magnitude starting from the zero vector;
\item the CPU time cost in seconds for completing the preorder phase (denoted by $t_{p}$), for constructing the approximate inverse factorization ($t_{f}$), and for solving the linear system ($t_{s}$). Symbol ``-'' means that the corresponding phase does not apply to the given run. For example, some of the preconditioners used for the comparison against our method do not have a preorder phase.
\end{enumerate}

The codes were developed in Fortran 95 and the experiments were run in double precision floating point arithmetic on a PC equipped with an Intel(R) Core(TM)2 Duo CPU E8400, 3.00 GHz of frequency, 4 GB of RAM and 6144 KB of cache memory. In the coming sections we study the effect of using different parameter settings, and we illustrate the overall performance on the selected matrix problems.

\medskip

\subsection{Performance of the multilevel mechanism}~\label{sec:5.1}

The AMES method can be seen as a multilevel generalization of factorized approximate inverse techniques such as the FSAI preconditioner by Kolotilina and Yeremin, and the AINV preconditioner 
by Benzi and T\accent23uma, that we recalled in Section~\ref{sec:3}. Therefore, first we present some comparison between these methods, to show the benefit of the multilevel mechanism. The results are reported in Table~\ref{tab:fsaiainv}. For these runs, we considered four matrix problems from Table~\ref{tab:problems}, that are \emph{orsirr\_1}, \emph{1138\_bus}, \emph{bcsstk27} and \emph{epb0}. In our AMES solver, we inverted the last level blocks using ILU, FSAI and AINV factorizations. For ILU, we used the multilevel implementation available in the ILUPACK package~\cite{ilupack:2010} ~(this combination is referred to as \textit{AMES\_ILU} in the table). For FSAI, we used the structure of the nonzero pattern of the lower (resp. upper) triangular part of the symmetrized block for the approximate inverse factors, and also the square of this pattern (this combination is referred to as \textit{AMES\_FSAI}). Finally, for AINV we used the implementation kindly provided by the authors (this combination is referred to as \textit{AMES\_AINV}). The dropping threshold value selected for the \textit{AINV}, \textit{AMES\_AINV} and \textit{AMES\_ILU} methods (referred to as {\it{Drop}} in the Table) is an absolute value, and was computed so that the resulting preconditioners had roughly equal memory cost. We used the default value for the parameter $condest=10$ (norm bound for the inverse factors $L^{-1}$ and $U^{-1}$) in ILUPACK.

\begin{figure}[ht!]
\begin{subfigmatrix}{3}
\subfigure[$L^{-1}$]%
  {\includegraphics[width=5cm]{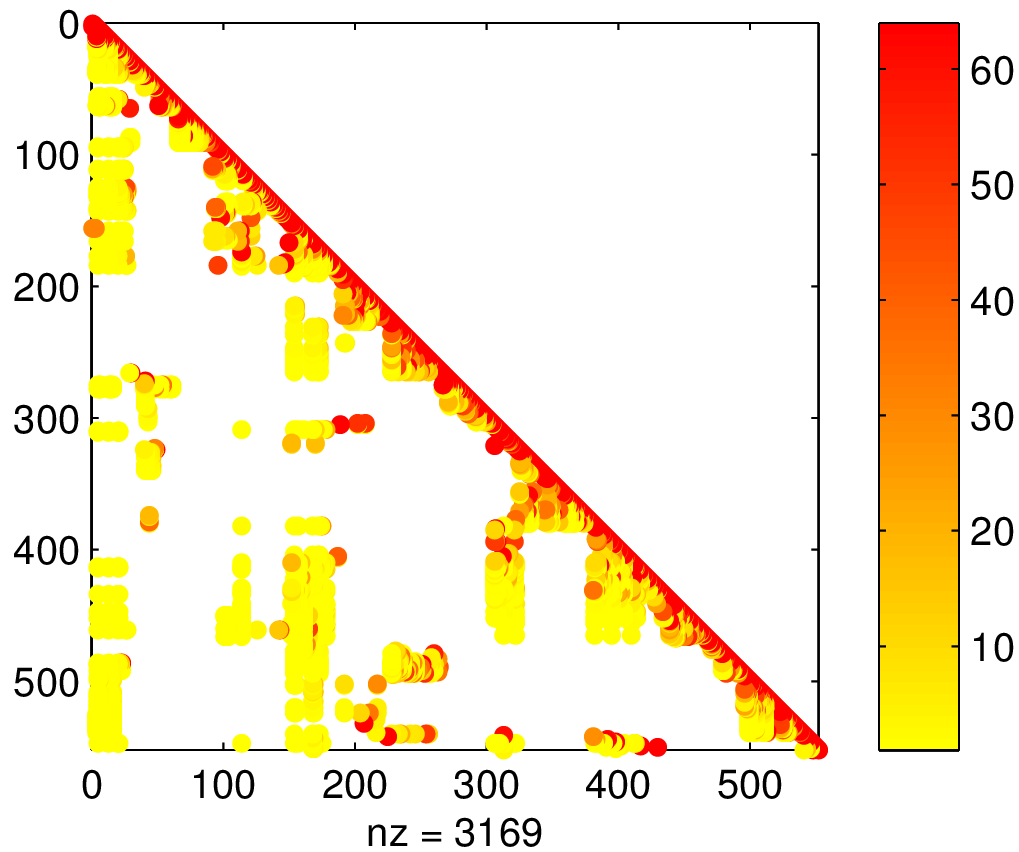}}
\subfigure[$M_L$]%
  {\includegraphics[width=5cm]{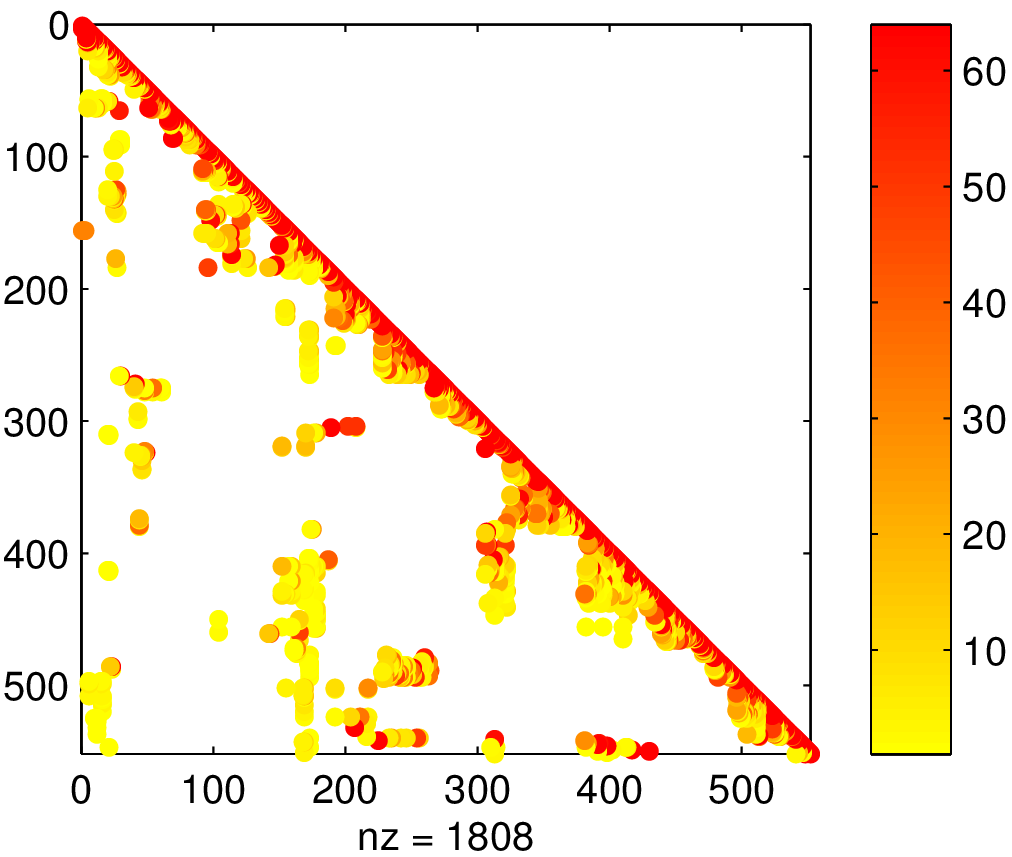}}
\subfigure[$W^T$]%
  {\includegraphics[width=5cm]{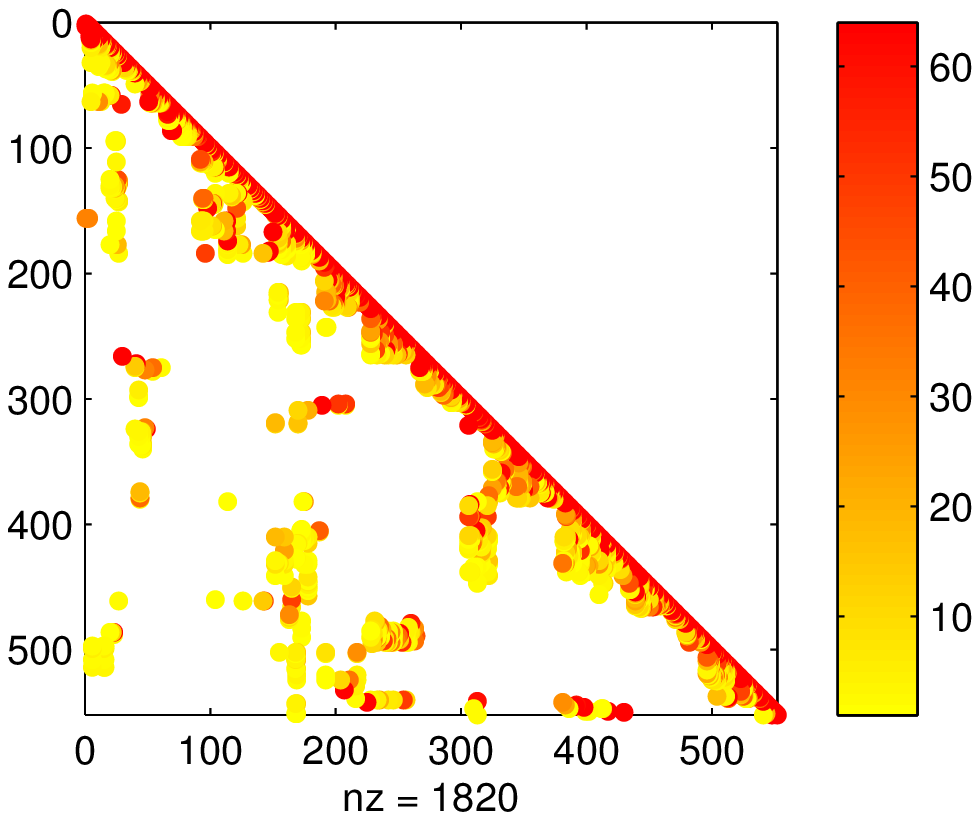}}
\end{subfigmatrix}
  \caption{The exact~(left) and approximate~(middle and right) inverse lower triangular factors of the \emph{1138\_bus} matrix.}~\label{fig:fsai1}
\end{figure}

\begin{figure}[ht!]
\begin{subfigmatrix}{3}
\subfigure[$L^{-1}$]%
  {\includegraphics[width=5cm]{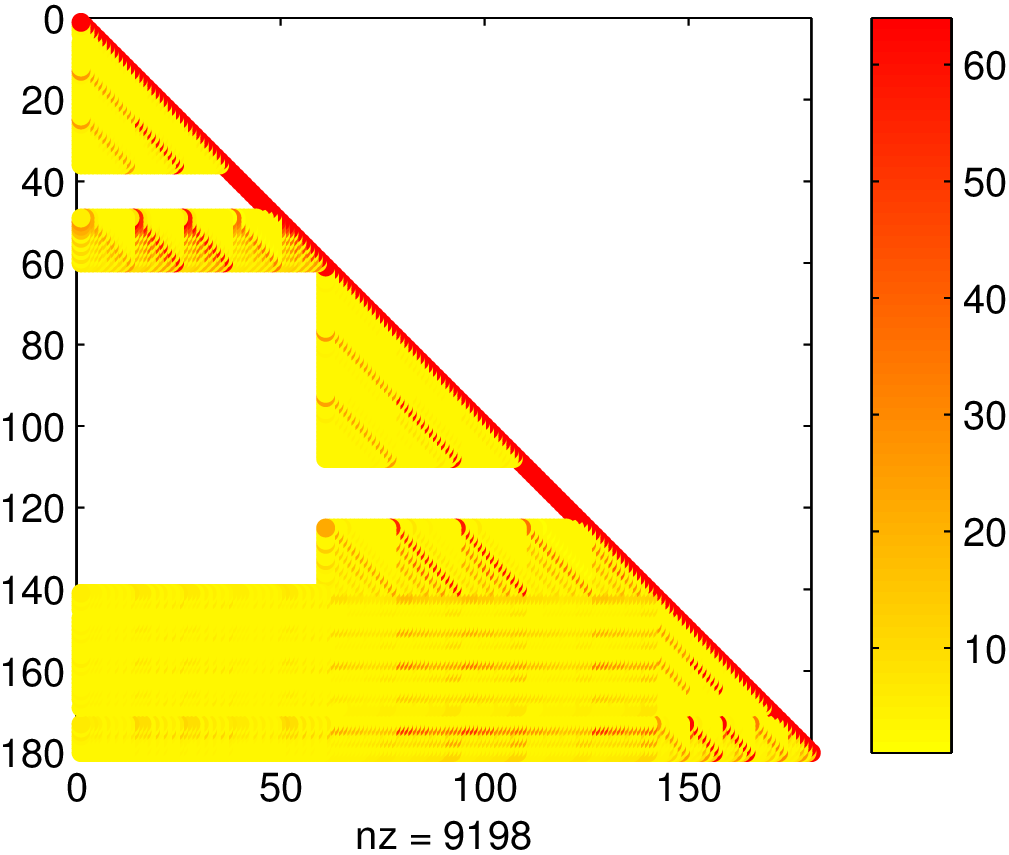}}
\subfigure[$M_L$]%
  {\includegraphics[width=5cm]{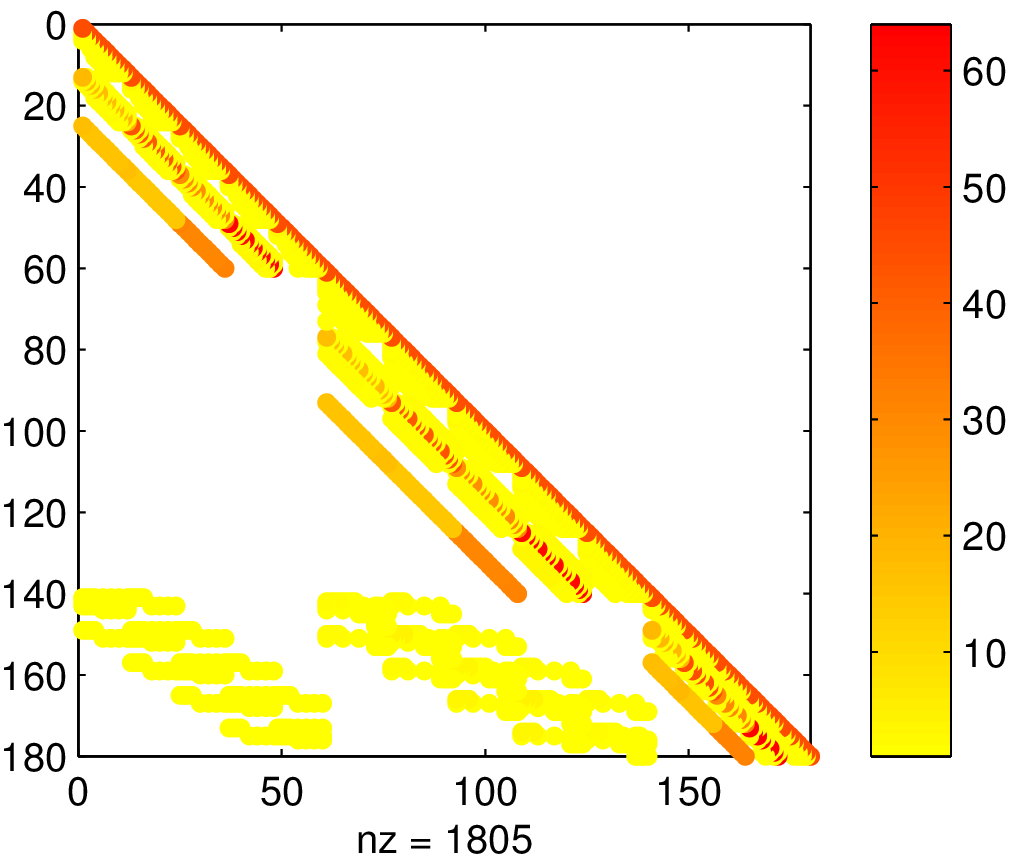}}
\subfigure[$W^T$]%
  {\includegraphics[width=5cm]{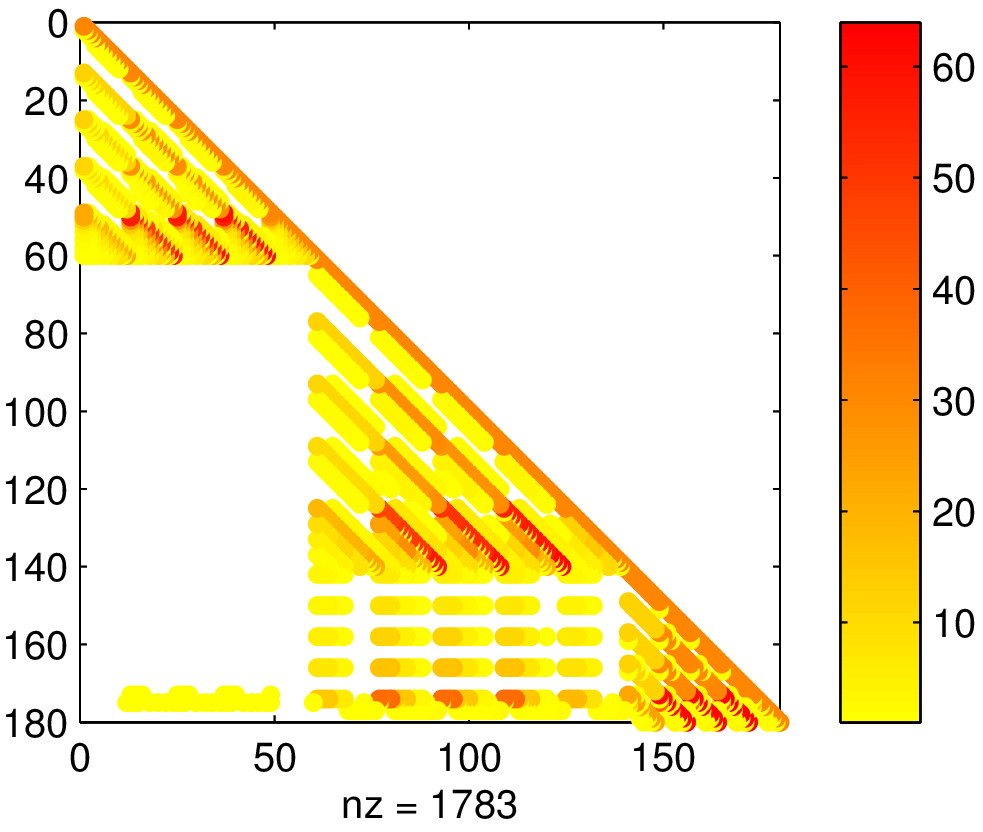}}
\end{subfigmatrix}
  \caption{The exact~(left) and approximate~(middle and right) inverse lower triangular factors of the \emph{orsirr\_1} matrix.}~\label{fig:fsai2}
\end{figure}

\begin{figure}[ht!]
\begin{subfigmatrix}{3}
\subfigure[$L^{-1}$]%
  {\includegraphics[width=5cm]{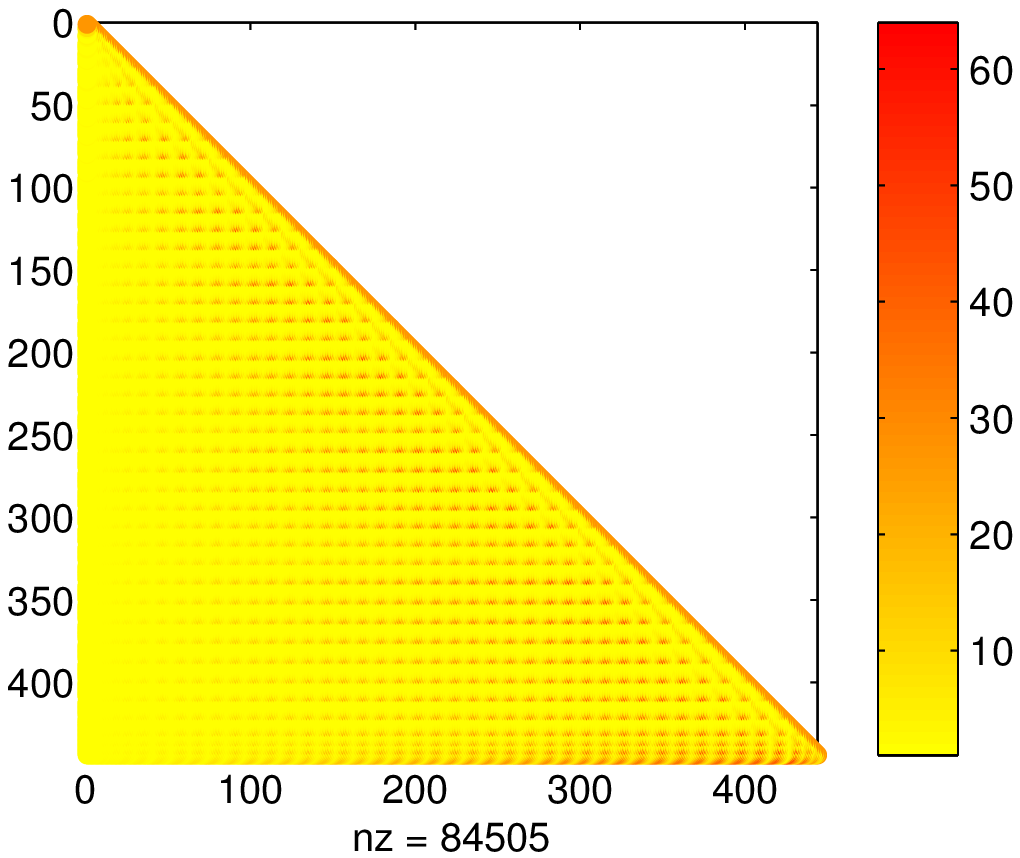}}
\subfigure[$M_L$]%
  {\includegraphics[width=5cm]{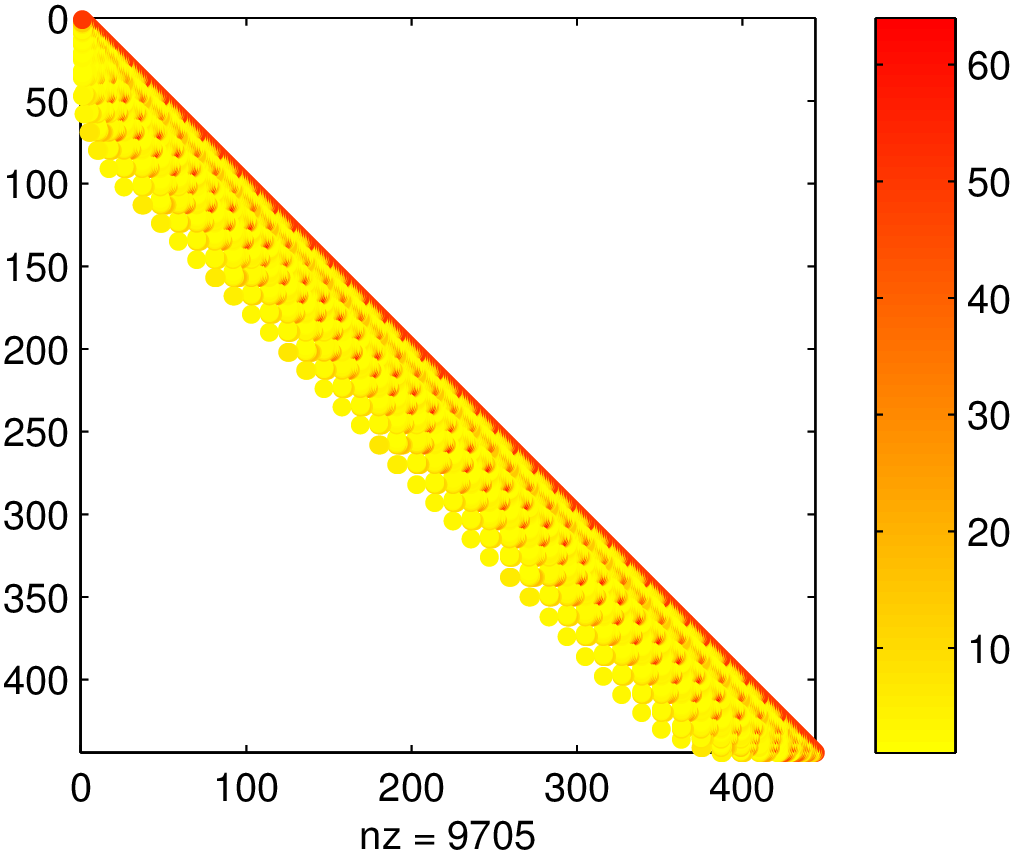}}
\subfigure[$W^T$]%
  {\includegraphics[width=5cm]{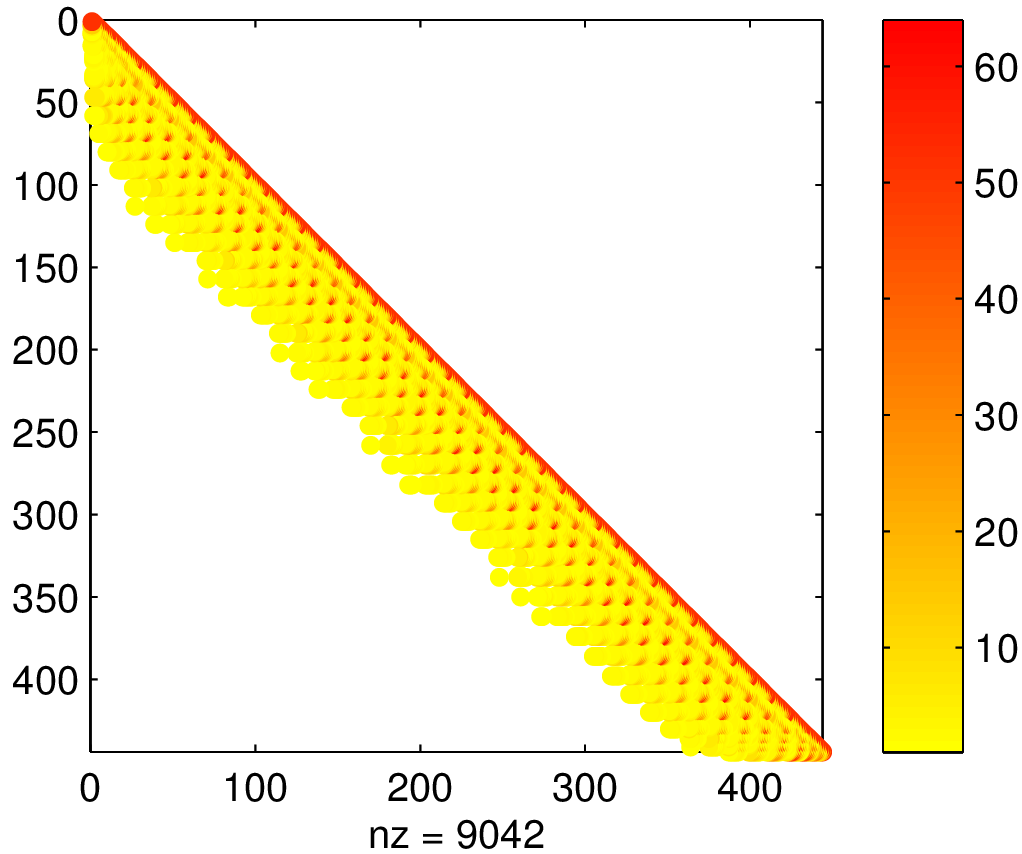}}
\end{subfigmatrix}
  \caption{The exact~(left) and approximate~(middle and right) inverse lower triangular factors of the \emph{cz40948} matrix.}~\label{fig:fsai3}
\end{figure}

\begin{figure}[ht!]
\begin{subfigmatrix}{3}
\subfigure[$L^{-1}$]%
  {\includegraphics[width=5cm]{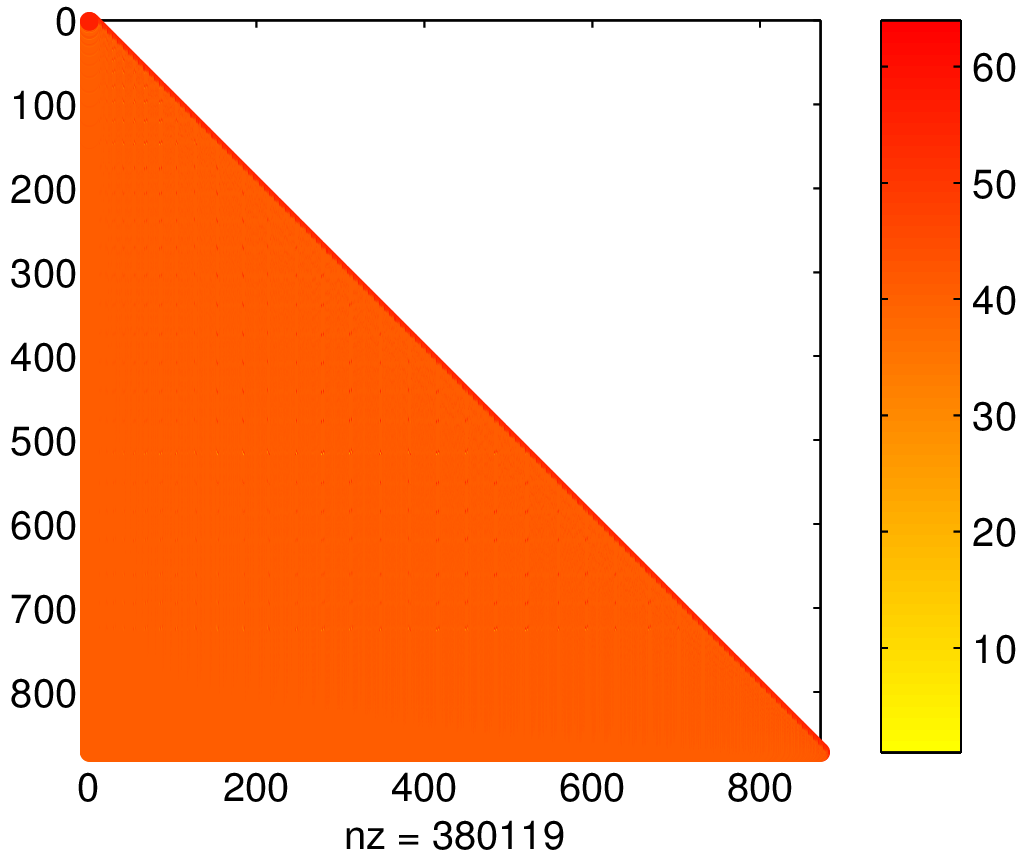}}
\subfigure[$M_L$]%
  {\includegraphics[width=5cm]{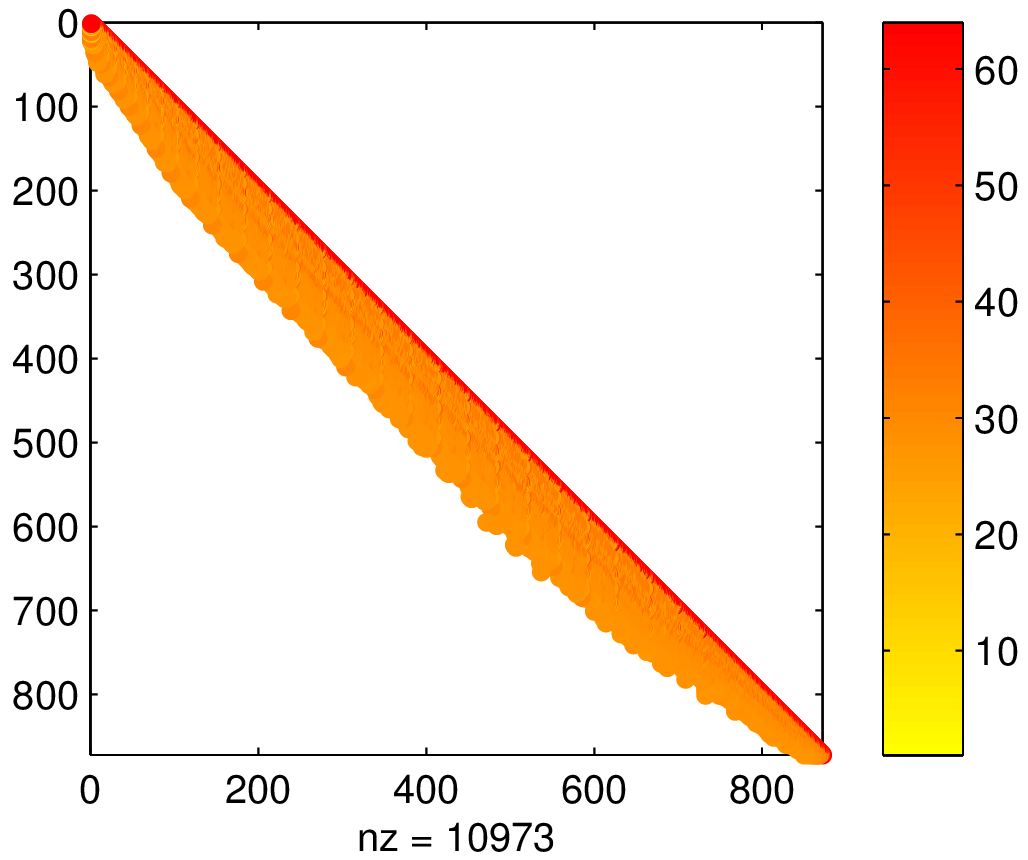}}
\subfigure[$W^T$]%
  {\includegraphics[width=5cm]{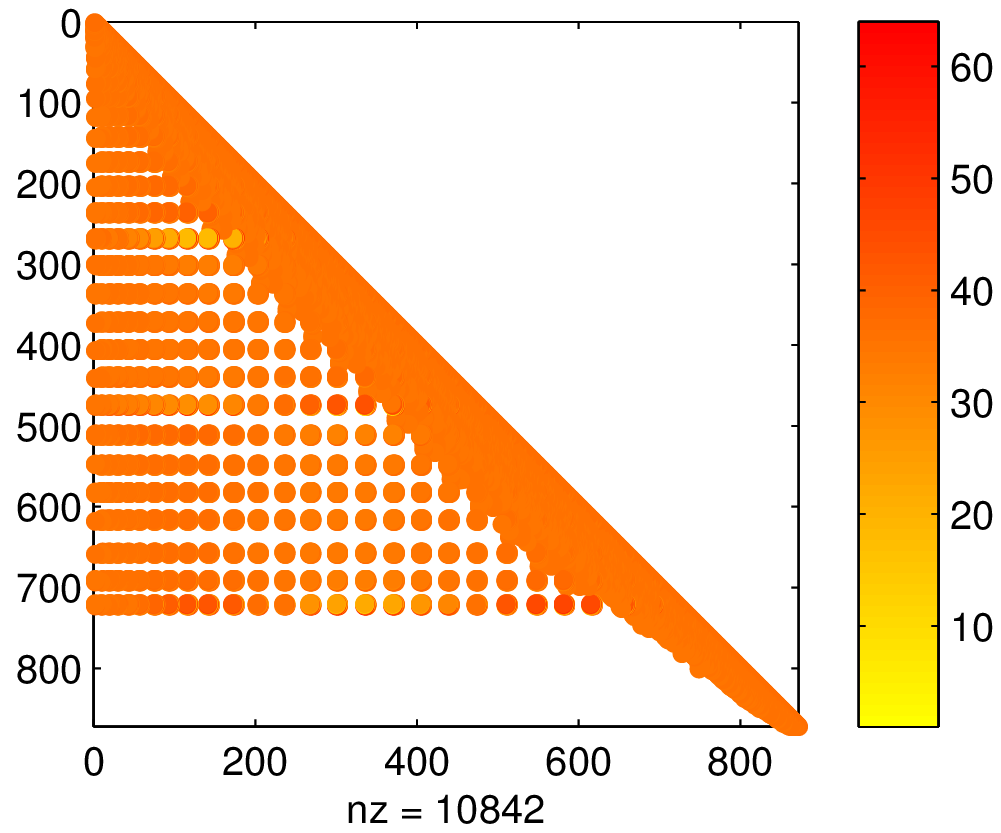}}
\end{subfigmatrix}
  \caption{The exact~(left) and approximate~(middle and right) inverse lower triangular factors of the \emph{ABACUS\_shell\_ud} matrix.}~\label{fig:fsai4}
\end{figure}

In our runs, the multilevel variants \textit{AMES\_FSAI} and \textit{AMES\_AINV} performed consistently better than the \textit{FSAI} and \textit{AINV} solvers in terms of convergence rate and storage cost. This is probably due to the multilevel mechanism that enabled us to exploit sparsity in the inverse factors more effectively. The best solutions with AMES were obtained using ILU as local solver, while the threshold-based dropping rules of the AINV method often computed a better pattern for the approximate inverse factors than the static approach used in the FSAI method. We can see evidence of this behaviour in Figures~\ref{fig:fsai1}~-~\ref{fig:fsai4}, where for one of the last-level blocks of the 
permuted coefficient matrix~(\ref{eq:perm}) we compare 
the structure of its exact inverse factor $L^{-1}$, and of the approximations 
$M_L$ and $W^T$ of  $L^{-1}$ as computed by, respectively, the \textit{AMES\_FSAI} code using the square of the pattern of the symmetrized block, and by the \textit{AMES\_AINV} code. Large to small entries are depicted in different colors, from red to orange and yellow. The approximation is good  for the \emph{1138\_bus} problem (Figure~\ref{fig:fsai1}) but poor for the \textit{orsirr\_1} matrix (Figure~\ref{fig:fsai2}), and this is confirmed by the different convergence results for the two problems, reported in Table~\ref{tab:fsaiainv}. 
On some larger problems, like the \emph{cz40948} and the \textit{ABACUS\_shell\_ud} problems, shown in Figures~\ref{fig:fsai3}~-~\ref{fig:fsai4}, we found that $L^{-1}$ had no evident structure; in this case we had to increase the number of nonzeros in $M_L$ and $W^T$ significantly to converge. For example on the \textit{ABACUS\_shell\_ud} problem, \textit{AMES\_AINV} converged in 468 iteration with $\frac{nnz(Z+W)}{nnz(A)}=11.6$ while \textit{AMES\_FSAI} did not converge at this value of density.
In these situations, uniformly better convergence were obtained using ILU as local solver. We will focus mostly on this choice of local solver for the experiments of this paper. Notice that in this case the entries of the inverse factors are not computed explicitly, and the application of the preconditioner is carried out through a backward and forward substituion procedure. 
Other options may be considered for the last level solver, such as the ARMS method~\cite{ARMS} and enhanced FSAI methods~\cite{EBFSAI}, but these are not included in the presented analysis.

\begin{table}[!ht]
\begin{center}

\subtable[orsirr\_1]{
   \begin{tabular}{lccccccc}
    \hline
 Method& ~~Pattern~~ & Drop/condest & ~\emph{Its}~ & $\frac{nnz(M_L+M_U)}{nnz(A)}$ & ~$t_p$~ & ~$t_f$~ & ~$t_s$~ \\

\hline
 \tabincell{l}{\textit{AMES\_FSAI}\\\textit{FSAI}} & $A+A^{T}$  & -&
\tabincell{c}{273 \\ 304} & \tabincell{c}{1.42 \\ 1.45 } & \tabincell{c}{0.011 \\-} & \tabincell{c}{0.023 \\ 0.070 } & \tabincell{c}{0.22 \\ 0.23 }\\

\hline
\tabincell{l}{\textit{AMES\_FSAI}\\\textit{FSAI}}  &  $(A+A^{T})^{2}$ &-&
\tabincell{c}{217 \\ 236} & \tabincell{c}{3.43 \\ 3.76} & \tabincell{c}{0.013 \\-} & \tabincell{c}{0.035 \\ 0.088 } & \tabincell{c}{0.17 \\ 0.16 }\\

\hline
\tabincell{l}{\textit{AMES\_AINV}\\\textit{AINV}}  & - & \tabincell{c}{0.03 \\ 0.07} &
\tabincell{c}{67 \\ 80} & \tabincell{c}{2.27 \\ 2.22 } & \tabincell{c}{0.016 \\-} & \tabincell{c}{0.014 \\ 0.016 } & \tabincell{c}{0.034 \\ 0.024 }\\

\hline
\textit{AMES\_ILU}  & - & \tabincell{c}{8e-3/10}  & \tabincell{c}{31} & \tabincell{c}{1.24}& \tabincell{c}{0.012} & \tabincell{c}{0.013} & \tabincell{c}{7.4e-3}\\
\hline
\end{tabular}
}

\subtable[\texttt{1138\_bus}]{
  \begin{tabular}{lccccccc}
    \hline
 Method& ~~Pattern~~ & Drop/condest & ~\emph{Its}~ & $\frac{nnz(M_L+M_U)}{nnz(A)}$  & ~$t_p$~ & ~$t_f$~ & ~$t_s$~ \\

\hline
 \tabincell{l}{\textit{AMES\_FSAI}\\\textit{FSAI}} & $A+A^{T}$  &-&
\tabincell{c}{7 \\ 9} & \tabincell{c}{2.24 \\ 2.32 } & \tabincell{c}{5.2e-3 \\-} & \tabincell{c}{0.032 \\ 0.074 } & \tabincell{c}{1.2e-3 \\ 8.9e-4 }\\

\hline
\tabincell{l}{\textit{AMES\_FSAI}\\\textit{FSAI}}  &  $(A+A^{T})^{2}$ &-&
\tabincell{c}{5 \\ 6} & \tabincell{c}{2.70 \\ 2.88 } & \tabincell{c}{5.0e-3 \\-} & \tabincell{c}{0.035 \\ 0.077 } & \tabincell{c}{1.0e-3 \\ 6.4e-4 }\\

\hline
\tabincell{l}{\textit{AMES\_AINV}\\\textit{AINV}}  & - &\tabincell{c}{0.6 \\ 0.7} &
\tabincell{c}{13 \\ 16} & \tabincell{c}{2.85 \\ 2.88 } & \tabincell{c}{7.0e-3 \\-} & \tabincell{c}{2.0e-3 \\ 6.0e-3 } & \tabincell{c}{1.9e-3 \\ 3.2e-3 }\\

\hline
\textit{AMES\_ILU}  & - & 0/10 & 1 & 1.00 & \tabincell{c}{5.1e-3} & \tabincell{c}{3.9e-3} & \tabincell{c}{7.0e-4}\\
\hline
\end{tabular}
}

\end{center}
\end{table}

\begin{table} 
\begin{center} 
\subtable[bcsstk27]{
   \begin{tabular}{lccccccc}
    \hline
 Method& ~~Pattern~~ & Drop/condest & ~\emph{Its}~ & $\frac{nnz(M_L+M_U)}{nnz(A)}$  & ~$t_p$~ & ~$t_f$~ & ~$t_s$ ~\\

\hline
 \tabincell{l}{\textit{AMES\_FSAI}\\\textit{FSAI}} & $A+A^{T}$  &-&
\tabincell{c}{8 \\ 19} & \tabincell{c}{0.90 \\ 1.27 } & \tabincell{c}{0.062 \\-} & \tabincell{c}{0.041 \\ 0.20 } & \tabincell{c}{0.021 \\ 4.1e-3 }\\

\hline
\tabincell{l}{\textit{AMES\_FSAI}\\\textit{FSAI}}  &  $(A+A^{T})^{2}$ &-&
\tabincell{c}{5 \\ 13} & \tabincell{c}{1.16 \\ 2.72 } & \tabincell{c}{0.063 \\-} & \tabincell{c}{0.071 \\ 0.47 } & \tabincell{c}{0.018 \\ 3.7e-3 }\\

\hline
\tabincell{l}{\textit{AMES\_AINV}\\\textit{AINV}}  & - &\tabincell{c}{1e-3 \\ 0.06} &
\tabincell{c}{6 \\ 16} & \tabincell{c}{1.18 \\ 0.98 } & \tabincell{c}{0.055 \\-} & \tabincell{c}{0.040 \\ 0.063 } & \tabincell{c}{7.3e-3 \\ 5.7e-3 }\\

\hline
\textit{AMES\_ILU}  & - & \tabincell{c}{0.01/10} & \tabincell{c}{6} & \tabincell{c}{0.978} & \tabincell{c}{0.059} & \tabincell{c}{0.016} & \tabincell{c}{0.010}\\
\hline
\end{tabular}
}

\subtable[epb0]{
   \begin{tabular}{lccccccc}
    \hline
 Method& ~~Pattern~~ & Drop/condest & ~\emph{Its}~ & $\frac{nnz(M_L+M_U)}{nnz(A)}$  & ~$t_p$~ & ~$t_f$~ & ~$t_s$~ \\

\hline
 \tabincell{l}{\textit{AMES\_FSAI}\\\textit{FSAI}} & $A+A^{T}$  &-&
\tabincell{c}{277 \\ 400} & \tabincell{c}{1.67 \\ 1.69 } & \tabincell{c}{0.020 \\-} & \tabincell{c}{0.011 \\ 0.19 } & \tabincell{c}{0.66 \\ 0.59 }\\

\hline
\tabincell{l}{\textit{AMES\_FSAI}\\\textit{FSAI}}  &  $(A+A^{T})^{2}$ &-&
\tabincell{c}{161 \\ 290} & \tabincell{c}{4.32 \\ 4.81} & \tabincell{c}{0.021 \\-} & \tabincell{c}{0.023 \\ 0.23 } & \tabincell{c}{0.40 \\ 0.27 }\\

\hline
\tabincell{l}{\textit{AMES\_AINV}\\\textit{AINV}}  & - &\tabincell{c}{0.5 \\ 0.9} &
\tabincell{c}{132 \\ 347} & \tabincell{c}{3.32 \\ 4.26 } & \tabincell{c}{0.024 \\-} & \tabincell{c}{4.5e-3 \\ 0.015 } & \tabincell{c}{0.21 \\ 0.42 }\\

\hline
\textit{AMES\_ILU}  & - & \tabincell{c}{0.1/10} & \tabincell{c}{7} & \tabincell{c}{1.848} & \tabincell{c}{0.020} & \tabincell{c}{4.1e-3} & \tabincell{c}{0.019}\\
\hline
\end{tabular}
}

\end{center}
\vspace{1cm}
\caption{Numerical experiments on selected matrix problems illustrating the performance of the multilevel sparse approximate inverse preconditioner AMES against the factorized approximate inverse methods FSAI and AINV.}\label{tab:fsaiainv}
\end{table}

\subsection{Varying the number of independent clusters at the first level}~~\label{sec:5.2}
We considered three matrix problems in our runs: \emph{cz20468}, \emph{ABACUS\_shell\_ud} and \emph{cz40948}. In Table~\ref{tab:number p} we show the results varying the number of independent clusters $p$ at the first level of reordering of $A$ in~(\ref{eq:perm}). For each problem, we used the same number of levels $n_{lev}$ in the AMES structure, and tuned the drop tolerance in the local ILU factorization to keep the memory ratio $\frac{nnz(M_L+M_U)}{nnz(A)}$ roughly constant while increasing $p$ in different runs.
Clearly, larger $p$ results in more independent clusters of smaller size, and in larger Schur complement matrices. In the table we report the ratio $\frac{sizeB}{sizeA_{S}}$ between the average size of the independent clusters $B_i$ and the size of the Schur complement at the first level. Increasing $p$ reduces in turn $\frac{sizeB}{sizeA_{S}}$ to values smaller than 1. Using ILU as local solver, the best convergence results were obtained when $\frac{sizeB}{sizeA_{S}} \approx 1$. 
Our experiments indicate that for good performance the size of each independent cluster should be approximately equal to that of the Schur complement. 

\begin{table}[!ht]
\begin{center}

\begin{tabular}{lccccccc}
\hline
Matrix &$~~p~~$ & ~\emph{Its}~  & ~~$t_p$~~ & ~$t_f$~  &~~$t_s$~~ & ~$t_{per\_it}$ ~& $~\frac{sizeB}{sizeA_{S}}$~\\
\hline

\texttt{cz20468}  & \tabincell{c}{ 15\\20\\30 \\ 40} &\tabincell{c}{ 151 \\ 139 \\ 131 \\ 137 } & \tabincell{c}{ 0.4\\0.4\\0.5\\0.5 } & \tabincell{c}{0.3\\0.3\\0.4\\0.6} &\tabincell{c}{ 3.0\\2.8\\2.5\\2.6} & \tabincell{c}{0.020\\0.020\\0.019\\0.019} & \tabincell{c}{4.3\\2.5\\1.1\\0.6}\\

\hline
\texttt{ABACUS\_shell\_ud}  & \tabincell{c}{ 4\\6\\12 \\ 15} &\tabincell{c}{ 258 \\ 242 \\213\\253 } & \tabincell{c}{ 0.5\\0.5\\0.5\\0.5 } & \tabincell{c}{1.2\\1.2\\2.0\\2.1 } &\tabincell{c}{ 9.1\\8.3\\6.8\\8.7} & \tabincell{c}{0.035\\0.034\\0.032\\0.034} & \tabincell{c}{7.8\\3.8\\1.0\\0.7}\\
\hline

\texttt{cz40948}  & \tabincell{c}{ 15\\30 \\45\\ 50} &\tabincell{c}{ 219\\212\\198\\207 } & \tabincell{c}{ 1.1\\1.1\\1.1\\1.1 } & \tabincell{c}{0.5\\0.7\\1.2\\1.9} &\tabincell{c}{ 10.8\\10.0\\9.2\\9.8} & \tabincell{c}{0.049\\0.047\\0.046\\0.047} & \tabincell{c}{8.7\\2.2\\0.9\\0.5}\\
\hline

\end{tabular}
\end{center}
\caption{The best performance of the multilevel sparse approximate inverse preconditioner are observed when $~\frac{sizeB}{sizeA_{S}}\approx 1$.}\label{tab:number p}

\end{table}

\subsection{Varying the number of reduction levels for the diagonal blocks}~\label{sec:5.3}
We consider again matrices \emph{cz40948}, \emph{ABACUS\_shell\_ud} and \emph{cz20468} for our tests. We varied the number of levels $n_{lev}$ from 1 to 3
in the multilevel reordering of the diagonal blocks. In each run we tuned
the dropping threshold parameter to have roughly the same memory cost in the experiments for each matrix. We chose the value of $p$ for each problem so that $\frac{sizeB}{sizeA_{S}}\approx1$ as reported in Section 5.2.
The last level blocks were factorized using ILUPACK~\cite{ilupack:2010}. The results reported in Table~\ref{tab:level1} show that using more levels can reduce the number of iterations for similar memory ratio as we can gain additional sparsity during the factorization. However, probably due to our non optimized implementation, the solution cost tends to increase. From our experiments, a small number of reduction levels is recommended to use.
\begin{table}[!ht]
\begin{center}
\begin{tabular}{lcccccc}
\hline
Matrix & ~$\emph{$n_{lev}$}$~ & ~\emph{Its}~  & ~$t_p$~ & ~$t_f$~ & ~$t_s$~ & ~$t_{per\_it}$~\\

\hline

\texttt{cz20468}  & \tabincell{c}{1\\ 2\\ 3}  &
\tabincell{c}{113\\80\\71} & \tabincell{c}{0.5\\0.5\\0.6 } & \tabincell{c}{0.4\\0.5\\0.5} &
\tabincell{c}{1.9\\2.3\\4.1}  & \tabincell{c}{0.017\\0.028\\0.058} \\

\hline

\texttt{ABACUS\_shell\_ud}  & \tabincell{c}{1\\ 2\\ 3}  &
\tabincell{c}{388\\381\\294} & \tabincell{c}{0.5 \\0.5 \\ 0.6 } & \tabincell{c}{1.7\\1.9\\1.9} &
\tabincell{c}{17.6\\21.9\\22.7}  & \tabincell{c}{0.045\\0.057\\0.077} \\

\hline

\texttt{cz40948}  & \tabincell{c}{1\\ 2\\ 3}  &
\tabincell{c}{198\\154\\133} & \tabincell{c}{1.1 \\1.2 \\ 1.3 } & \tabincell{c}{1.2\\1.3\\1.3} &
\tabincell{c}{9.2\\10.4\\17.3 }  & \tabincell{c}{0.046\\0.068\\0.13} \\

\hline

\end{tabular}
\vspace{1cm}
\caption{The number of iterations of the multilevel approximate inverse preconditioner can be reduced by increasing the number of reduction levels $n_{lev}$ for the diagonal blocks, at roughly equal memory costs.}\label{tab:level1}
\end{center}
\end{table}

\subsection{Varying the number of reduction levels for the Schur complement}~~\label{sec:5.4}
The Schur complement matrix relative to the block $C$ in~(\ref{eq:perm}) typically preserves a good deal of sparsity, and this can be further exploited during the factorization by applying, e.g., the multilevel nested dissection reordering to $A_S$, similarly to what is done to the upper leftmost block $B$. We implemented this idea at the first permutation level, using ILU factorization as local solver and selecting the same values of $p$ and $n_{lev}$ for each matrix problem. We tuned the drop tolerence in the ILU factorization to have roughly the same memory costs in different runs. We varied $n_{levAS}$ from 0 to 3 ($n_{levAS}=0$ means that only the diagonal blocks of the upper-left block B are permuted). The results reported in Table~\ref{tab:reduction_schur} show that the simultaneous permutation of both the diagonal blocks of $B$ and of the Schur complement can make the preconditioner more robust. We adopted this strategy in the experiments illustrated in the coming sections, selecting in each run the value of $n_{levAS}$ that minimized the total solution cost.

\begin{table}[!ht]
\begin{center}
\begin{tabular}{lcccccc}
\hline
Matrix & ~$\emph{$n_{levAS}$}$~ & ~\emph{Its}~  & ~$t_p$~ & ~$t_f$~ & ~$t_s$~ \\

\hline

\texttt{cz20468}  & \tabincell{c}{0\\1\\ 2\\ 3}  &
\tabincell{c}{331\\228 \\ 209 \\181} & \tabincell{c}{0.5\\0.4 \\ 0.4  \\ 0.4 } & \tabincell{c}{0.2\\1.3 \\ 1.3 \\ 1.3} & \tabincell{c}{8.2\\5.6 \\ 4.8\\4.0 } \\

\hline

\texttt{ABACUS\_shell\_ud}  & \tabincell{c}{0\\1\\ 2\\ 3} &
\tabincell{r}{576 \\485 \\ 414 \\393} & \tabincell{c}{0.5 \\ 0.5 \\ 0.5  \\ 0.5 } & \tabincell{r}{1.8 \\1.8 \\ 1.4 \\ 1.6} & \tabincell{r}{35.0 \\29.5 \\ 24.4\\22.2 } \\
\hline

\texttt{cz40948}  & \tabincell{c}{0\\1\\ 2\\ 3}  &
\tabincell{c}{183\\166 \\ 157 \\ 152} & \tabincell{c}{1.9 \\1.9 \\ 1.9  \\ 1.8 } & \tabincell{c}{0.5 \\6.4  \\ 6.1 \\ 6.1} & \tabincell{c}{23.7 \\16.6 \\ 14.8 \\ 14.3 } \\

\hline

\end{tabular}
\end{center}
\vspace{1cm}
\caption{At roughly equal memory costs, larger reduction levels for the Schur complement can improve the convergence rate.}\label{tab:reduction_schur}
\end{table}

\subsection{Comparison against other solvers}~
We compared the performance of the AMES preconditioner against other three popular algebraic preconditioners for solving linear systems, that are the ILUPACK solver developed by Bollh{\"o}fer and Saad~\cite{ilupack:2010}, the Algebraic Recursive Multilevel Method (ARMS) proposed by Saad and Suchomel~\cite{ARMS}, and the SParse Approximate Inverse preconditioner (SPAI) introduced by Grote and Huckle~\cite{grhu:97}. 
As in the previous experiments, for each run we recorded the CPU time from the start of the solution until the initial residual was reduced by 12 orders of magnitude or until the process failed. We declared a solver failure when no convergence was achieved after 5000 iterations of the restarted GMRES method. We selected the parameters carefully to have a fair comparison between different methods. In AMES, following our conclusions from Section~\ref{sec:5.2}, we selected the number of blocks $B_i$ at the first level so that their average size is almost equal to the size of the Schur complement. For every problem we tested different combinations of number of levels $n_{lev}$ of recursive factorization and different values for the dropping threshold parameter $droptol$ for the factorization of the last level blocks $B_i$ and of the Schur complements. We chose the best combination in terms of memory and time to solution costs for the given problem. Then we tuned the value of the dropping threshold in the ILUPACK, ARMS, SPAI and AINV solvers to have roughly equal memory costs as in AMES, setting the other parameters equal to their default values defined in those packages. The performance of these methods is rather sensitive to the dropping threshold parameter. For example, on the $rma10$ problem, ILUPACK converged in only 9 iterations using the default value $droptol=0.01$, but the computation costed $\frac{nnz(M)}{nnz(A)}=8.9$ and $t_f=45s$; ARMS converged in 26 iterations with the default $droptol=0.001$, costing $\frac{nnz(M)}{nnz(A)}=33.9$ and $t_f=1111s$;
and SPAI could not converge in 5000 iterations with $\frac{nnz(M)}{nnz(A)}=0.19$, using the default value $droptol=0.6$.
The number of levels of recursive factorization in the multilevel methods ILUPACK and ARMS are calculated automatically by the original codes developed by their authors. 
We point out that the performance comparison between AMES and the other solvers at fixed memory occupation may be a little penalizing for the AINV, FSAI and SPAI preconditioners as
one-level approximate inverses inherently need more memory; the ARMS method is a multilevel solver, but it factorizes the diagonal blocks without any permutation.

In Table~\ref{tab:result1}, we show the complete results of our experiments. These include number of iterations (\emph{Its}), density ratio ($\frac{nnz(M_L+M_U)}{nnz(A)}$), time costs for the preordering ($t_p$), factorization~($t_f$) and solve phase~($t_s$). We also tested the unpreconditioned GMRES for these matrices problems, and no convergence is achieved. 
We clearly see the good potential of the multilevel mechanism incorporated in the AMES preconditioner to reduce the number of iterations of Krylov methods, also in comparison to other multilevel solvers at low to moderate memory costs.
In our examples, AMES was more robust than these solvers especially at low memory ratios.

\begin{table}[!ht]
\begin{center}

\subtable[cz20468]{
\begin{tabular}{lccrrr}
\hline
Method & $\frac{nnz(M_L+M_U)}{nnz(A)}$ & \emph{Its}  & $t_p$ & $t_f$ & $t_s$ \\
\hline
\tabincell{l}{AMES\\ ILUPACK\\ ARMS \\ SPAI } & \tabincell{c}{1.26 \\ 1.24 \\ 1.16 \\ 1.64} &
\tabincell{r}{187~ \\ 2500~ \\+5000~ \\ +5000~~} & \tabincell{c}{0.3 \\ - \\ - \\ -} & \tabincell{r}{0.2 \\ 0.4 \\ 0.1 \\ 4.0} &
\tabincell{r}{4.2 \\ 40.3 \\ +6.5 \\ +8.0}
\end{tabular}}

\subtable[raefsky3]{
\begin{tabular}{lccrrr}
\hline
Method & $\frac{nnz(M_L+M_U)}{nnz(A)}$ & \emph{Its}  & $t_p$ & $t_f$ & $t_s$ \\
\hline
\tabincell{l}{AMES\\ ILUPACK\\ ARMS \\ SPAI} & \tabincell{c}{0.54 \\ 0.55 \\ 2.38 \\ 1.83 } &
\tabincell{r}{235~ \\ 1224~ \\+5000~ \\ +5000~~} & \tabincell{c}{2.4 \\ - \\ - \\ -} & \tabincell{r}{3.7 \\ 2.8 \\ 2.4 \\ 5040} &
\tabincell{r}{10.0 \\ 25.2 \\ +23.5 \\ +243.0}
\end{tabular}}

\subtable[ABACUS\_shell\_ud]{
\begin{tabular}{lcccrrr}
\hline
Method & $\frac{nnz(M_L+M_U)}{nnz(A)}$ & \emph{Its}  & $t_p$ & $t_f$ & $t_s$ \\
\hline
\tabincell{l}{AMES\\ ILUPACK\\ ARMS \\  SPAI} & \tabincell{c}{1.79 \\ 1.82 \\ 1.88 \\ 2.41} &
\tabincell{r}{453~ \\ 1411~ \\+5000~ \\ +5000~~} & \tabincell{c}{0.3 \\ - \\ - \\ -} & \tabincell{r}{0.8 \\ 0.5 \\ 0.2 \\ 11.0} &
\tabincell{r}{22.1 \\ 26.6 \\ +7.6 \\ +12.0}
\end{tabular}}

\subtable[sme3Db]{
\begin{tabular}{lccrrr}
\hline
Method & $\frac{nnz(M_L+M_U)}{nnz(A)}$ & \emph{Its}  & $t_p$ & $t_f$ & $t_s$ \\
\hline
\tabincell{l}{AMES\\ ILUPACK\\ ARMS\\ SPAI} & \tabincell{c}{0.85 \\ 0.74 \\ 5.61 \\ 1.23} &
\tabincell{r}{407~ \\ 1210~ \\ +5000~ \\ +5000~~} & \tabincell{c}{3.5 \\ - \\ - \\ -} & \tabincell{r}{8.4 \\ 4.1 \\ 39.0 \\  3360 } &
\tabincell{r}{39.3 \\ 41.4 \\ +54.9 \\ +123.0}
\end{tabular}}
\end{center}
\end{table}

\begin{table} 
\begin{center} 
\subtable[viscoplastic2]{
\begin{tabular}{lccrrr}
\hline
Method & $\frac{nnz(M_L+M_U)}{nnz(A)}$ & \emph{Its}  & $t_p$ & $t_f$ & $t_s$ \\
\hline
\tabincell{l}{AMES\\ ILUPACK\\ ARMS \\  SPAI} & \tabincell{c}{3.07 \\ 4.00 \\ 3.02 \\ 3.37} &
\tabincell{r}{78~ \\ 2500~ \\ +5000~ \\ +5000~~} & \tabincell{c}{0.9 \\ - \\ - \\ -} & \tabincell{r}{14.3 \\ 1.6 \\ 0.9 \\  244.0} &
\tabincell{r}{3.9 \\ 70.0 \\ +10.9 \\ +24.0 }
\end{tabular}}

\subtable[cz40948]{
\begin{tabular}{lccrrr}
\hline
Method & $\frac{nnz(M_L+M_U)}{nnz(A)}$ & \emph{Its}  & $t_p$ & $t_f$ & $t_s$ \\
\hline
\tabincell{l}{AMES\\ ILUPACK\\ ARMS \\  SPAI} & \tabincell{c}{1.41 \\ 1.48 \\ 1.70 \\ 1.64} &
\tabincell{r}{170~ \\ 1627~ \\ +5000~ \\ +5000~~} & \tabincell{c}{0.7 \\ - \\ -  \\ -} & \tabincell{r}{0.4 \\ 1.0 \\ 0.9 \\ 8.5} &
\tabincell{r}{7.4 \\ 51.1 \\ +21.8 \\ +17.2}
\end{tabular}}

\subtable[rma10]{
\begin{tabular}{lcccrrr}
\hline
Method & $\frac{nnz(M_L+M_U)}{nnz(A)}$  & \emph{Its}  & $t_p$ & $t_f$ & $t_s$ \\
\hline
\tabincell{l}{AMES\\ ILUPACK\\ ARMS\\ SPAI} & \tabincell{c}{2.33 \\ 2.27 \\ 14.30 \\ 4.84} &
\tabincell{r}{164~ \\ ~~1242~ \\+5000~ \\ +5000~ } & \tabincell{c}{3.9 \\ - \\ - \\ -} & \tabincell{r}{13.1 \\ 8.6 \\ 203.9 \\ 11280} &
\tabincell{r}{34.5 \\ 82.9 \\ +111.3 \\ +180}
\end{tabular}}

\subtable[finan512]{
\begin{tabular}{lcccrrr}
\hline
Method & $\frac{nnz(M_L+M_U)}{nnz(A)}$ & ~\emph{Its}~  & ~$t_p$~ & ~$t_f$~ & ~$t_s$~ \\
\hline
\tabincell{l}{AMES\\ ILUPACK\\ ARMS\\ SPAI} & \tabincell{c}{0.59 \\ 0.62 \\ 0.58 \\ 0.61} &
\tabincell{r}{9~ \\ ~~11~ \\36~ \\ 7~ } & \tabincell{c}{0.8 \\ - \\ - \\ -} & \tabincell{r}{0.5 \\ 0.7 \\ 0.4 \\ 4.2} &
\tabincell{r}{0.8 \\ 0.1 \\ 0.5 \\ 0.2}
\end{tabular}}
\end{center}
\end{table}

\begin{table} 
\begin{center} 
\subtable[helm2d03]{
\begin{tabular}{lcccrrr}
\hline
Method & $\frac{nnz(M_L+M_U)}{nnz(A)}$& ~\emph{Its}~  & ~$t_p$~ & ~$t_f$~ & ~$t_s$~ \\
\hline
\tabincell{l}{AMES\\ ILUPACK\\ ARMS\\ SPAI} & \tabincell{c}{0.88 \\ 0.91 \\ 0.93 \\ 0.87} &
\tabincell{r}{6~ \\ ~~7~ \\12~ \\15~ } & \tabincell{c}{6.1 \\ - \\ - \\ -} & \tabincell{r}{4.3 \\ 3.7 \\ 1.4 \\ 100.7} &
\tabincell{r}{4.6 \\ 0.4\\ 1.5 \\2.7}
\end{tabular}}

\subtable[parabolic\_fem]{
\begin{tabular}{lcccrrr}
\hline
Method & $\frac{nnz(M_L+M_U)}{nnz(A)}$ & ~\emph{Its}~  & ~$t_p$~ & ~$t_f$~ & ~$t_s$~ \\
\hline
\tabincell{l}{AMES\\ ILUPACK\\ ARMS\\ SPAI} & \tabincell{c}{0.75 \\ 0.68 \\ 0.76 \\ 0.77} &
\tabincell{r}{4~ \\ 10~ \\~~12~ \\4~ } & \tabincell{c}{4.7 \\ - \\ - \\ -} & \tabincell{r}{5.7 \\ 5.3 \\ 2.0 \\ 175.3} &
\tabincell{r}{1.3 \\ 0.5 \\ 2.0 \\ 0.8}
\end{tabular}}

\end{center}
\vspace{1cm}
\caption{Performance comparison of the multilevel approximate inverse preconditioner against other iterative solvers, both one-level and multilevel.}\label{tab:result1}
\end{table}

\subsection{Effect of overlapping}~

We solved several problems from Table~\ref{tab:problems} combining the AMES method
with overlapping after the first level of reordering in~(\ref{eq:perm}).
In these runs, we set $n_{lev}=2$, and we tuned the $droptol$ parameter to have roughly the same memory costs in the experiments with and witout overlapping. In the last two columns of Table~\ref{tab:overlapping} we give the effect of overlapping on the change in size and in number of nonzeros for the overlapped system. The number of iteration ($\emph{Its}$) are almost the same after overlapping for problems $\emph{cz20468}$ and $\emph{cz40948}$, while for problems $\emph{sme3Db}$, $\emph{ABACUS\_shell\_ud}$ and $\emph{raefsky3}$ we observed a consistent reduction of the number of iterations $\emph{Its}$ by a factor between 9.5\% and 23.8\% and of the solving time $t_s$ by a factor between 21.4\% and 29.9\%. This is in agreement with our analysis of Section~\ref{sec:4}. In Table~\ref{tab:matrix_overlapping}, for each problem we studied the sparsity pattern of block $F$ and the size of blocks $B$ and $C$ before and after overlapping is applied at the first reordering level. The quantity $Sp_F$ denotes the ratio between the number of nonzero elements and the size of $F$, that is the sparsity degree $\frac{nnz(F)}{size(F)}$. As we can see, the $\emph{cz20468}$ and $\emph{cz40948}$ problems have the smallest relative size of the separator $C$ and also the smallest value of $Sp_{F}$; this means that less information is added to the subdomains. Following the analysis reported in Section~\ref{sec:4}, the overlapping technique is less likely to help on these two matrices, and this is also confirmed by the numerical results. Differently, problems $\emph{sme3Db}$ and $\emph{raefsky3}$ show larger values of $size_C$ and $Sp_F$ and in fact overlapping has a better effect on convergence for these two problems. In our experiments we found that a small number of independent clusters $p$ is recommended to use when overlapping.

\begin{table}[!ht]
\begin{center}
\begin{tabular}{lccrrcc}
\hline
Matrix & Method & $\frac{nnz(M_L+M_U)}{nnz(A)}$ & \emph{Its}  & $t_s$ &$\frac{n(A_{overlapped})}{n(A)}$ & $\frac{nnz(A_{overlapped})}{nnz(A)}$ \\

\hline
\texttt{cz20468}  & \tabincell{l}{overlapping\\ without overlapping} & \tabincell{c}{1.33\\ 1.32} &
\tabincell{r}{147\\ 149} & \tabincell{l}{4.3\\4.5 } &\tabincell{c}{1.005 \\ -} &\tabincell{c}{1.004 \\ -}\\

\hline
\texttt{raefsky3}  & \tabincell{l}{overlapping\\ without overlapping} & \tabincell{r}{0.56\\ 0.56} &
\tabincell{r}{218\\ 286} & \tabincell{r}{14.3\\ 20.3 } &\tabincell{c}{1.134 \\ -} &\tabincell{c}{1.135 \\ -}\\

\hline
\texttt{ABACUS\_shell\_ud}  & \tabincell{l}{overlapping\\ without overlapping} & \tabincell{r}{3.06\\ 3.03} &
\tabincell{r}{238\\ 263} & \tabincell{r}{13.6\\ 17.4 } &\tabincell{c}{1.020 \\ -} &\tabincell{c}{1.019 \\ -}\\

\hline
\texttt{sme3Db}  & \tabincell{l}{overlapping\\ without overlapping} & \tabincell{r}{0.91\\ 0.91} &
\tabincell{r}{389\\ 495} & \tabincell{r}{49.1\\62.5 } &\tabincell{c}{1.588 \\ -} &\tabincell{c}{1.639 \\ -}\\

\hline
\texttt{cz40948}  & \tabincell{l}{overlapping\\ without overlapping} & \tabincell{r}{1.42\\ 1.40} &
\tabincell{r}{175\\ 172} & \tabincell{r}{12.0\\17.5 } &\tabincell{c}{1.006 \\ -} &\tabincell{c}{1.004 \\ -}\\

\hline
\end{tabular}
\end{center}
\vspace{1cm}
\caption{Experiments on the effect of block overlapping on the performance of the multilevel sparse approximate inverse.}\label{tab:overlapping}
\end{table}

\begin{table}[!ht]
 \begin{center}
  \begin{tabular}{lccccrr}
    \hline
      Matrix problem~~ & Method & $size_B$&$size_C$~ & $Sp_F$~ \\
    \hline
     \texttt{cz20468} & \tabincell{l}{original\\ after overlapping} &\tabincell{c}{20405\\20450}  ~        &\tabincell{c}{63\\116}  ~    &\tabincell{c}{$1.1e-4$\\$4.3e-5$} \\
     
     \hline
     \texttt{raefsky3} & \tabincell{l}{original\\ after overlapping} &\tabincell{c}{19776\\21184}  ~        &\tabincell{c}{1424\\2864}  ~    &\tabincell{c}{$1.1e-4$\\$5.1e-4$}  \\
     
     \hline
     \texttt{ABACUS\_shell\_ud} & \tabincell{l}{original\\ after overlapping} &\tabincell{c}{23184\\23412}  ~        &\tabincell{c}{228\\458}  ~    &\tabincell{c}{$1.3e-4$\\$6.1e-5$}  \\
     
     \hline
     \texttt{sme3Db} & \tabincell{l}{original\\ after overlapping} &\tabincell{c}{19956\\25932}  ~        &\tabincell{c}{9111\\20214}  ~    &\tabincell{c}{$9.2e-4$\\$3.4e-4$}   \\
     
     \hline
     \texttt{cz40948} & \tabincell{l}{original\\ after overlapping}&\tabincell{c}{40825\\40925}  ~        &\tabincell{c}{123\\250}  ~    &\tabincell{c}{$5.2e-5$\\$2.4e-5$}    \\

    \hline
  \end{tabular}
  \ms
  \caption{\label{tab:matrix_overlapping} Effect of overlapping on matrix blocks size and sparsity.}
 \end{center}
\end{table}

\subsection{Utilizing direct solvers in the AMES framework}~

The results of previous sections indicate that the multilevel mechanism can be effective to reduce the memory burden but, at least in our implementation, tends to increase the cost per iteration. As an attempt of a possible remedy, we performed some runs setting the dropping threshold parameter $droptol$ equal to zero, and using a sparse direct solver, namely the routine {\texttt{MA38}} from the HSL Mathematical Software Library~\cite{hsl:2000}, as a local solver. No approximation is introduced and the Schur complements are exact. Therefore in each problem we can obtain convergence in one or two iterations, and the solving phase is much cheaper. This can be observed in Table~{\ref{tab:direct}} on selected matrix problems. Comparing against the results with inexact inversion, we see that using a direct solver as local component can save  computational time at only moderate extra storage cost.

\begin{table}[!ht]
 \begin{center}
\begin{tabular}{lcrrcc}
\hline
Matrix  & $\frac{nnz(M_L+M_U)}{nnz(A)}$ & \emph{Its}  & ~$t_p$~ & ~$t_f$~ & ~$t_s$~\\

\hline
\texttt{cz20468}   & \tabincell{c}{1.28} &
\tabincell{r}{2} & \tabincell{l}{0.7 } &\tabincell{c}{0.4} &\tabincell{c}{1.5}\\

\hline
\texttt{raefsky3}  & \tabincell{r}{2.74} &
\tabincell{r}{1} & \tabincell{r}{3.4 } &\tabincell{c}{11.1} &\tabincell{c}{1.3}\\

\hline
\texttt{cz40948} & \tabincell{c}{1.87} &
\tabincell{r}{2} & \tabincell{l}{1.2 } &\tabincell{c}{0.3} &\tabincell{c}{0.2}\\

\hline
\texttt{rma10}   & \tabincell{c}{3.01} &
\tabincell{r}{1} & \tabincell{l}{5.2 } &\tabincell{c}{11.6} &\tabincell{c}{0.8}\\

\hline
\end{tabular}
  \ms
  \caption{\label{tab:direct} Using the AMES factorization as a direct solver.}
 \end{center}
\end{table}

\section{Conclusions}\label{sec:6}

In this paper a recursive multilevel implementation of factorized sparse approximate inverse preconditioners for Krylov subspace methods was discussed. We used recursive combinatorial techniques and overlapping strategies as an attempt to remedy two typical drawbacks of explicit preconditioning, that are lack of robustness and high construction cost. 
The numerical experiments show that these strategies can improve the performance of conventional approximate inverse methods, yielding iterative solutions that can compete favourably against other popular solvers in use today. Parallelism can be exploited at various levels in our method, alongside other code optimization. Fine-grained blocking, filtering, postfiltering, adaptive pattern selection strategies have been shown to be promising approaches in other contexts~\cite{BFSAI,gBFSAI,JANNA,chow:01,chow:00a,VBARMS}, and these can be considered also in our setting. In a distributed memory implementation, it will be natural to split the oct-tree by assigning the local problems to different processors. 
An efficient use of recursive combinatorial algorithms may reduce considerably
the size of the Schur complements, hence the amount of inter-node communications. Memory demands, an important bottleneck of modern algorithms, are also limited, but this does not  penalize much the overall numerical efficiency of the solver, as illustrated by the experiments of Tables~\ref{tab:result1} and~\ref{tab:direct}. Overlapping does not destroy the sparsity structure of the matrix and can reduce further the interconnections between subdomains and separator set. Hence it is worthwhile considering it in a parallel setting as well. However, the parallel implementation of a fully distributed Schur complement formulation may not be trivial and will be considered in a separate study.

\section{Acknowledgements}~

The work of the first and of the third authors is funded by the Ubbo Emmius scholarship of the University of Groningen. The work of the last author is supported by NSFC (61170309), the Fundamental Research Funds for the Central Universities (ZYGX2013Z005). We are grateful to Miroslav T\r{u}ma for motivating discussion and for supplying the AINV package for our numerical experiments. Finally, we would like to thank the anonymous reviewers for their valuable and thorough comments that we believe helped improve significantly the presentation and the quality of the paper.

\end{document}